# EXPONENTIAL PENALTY FUNCTION CONTROL OF LOSS NETWORKS

By Garud Iyengar[1] and Karl Sigman[2]

*Columbia University*

We introduce penalty-function-based admission control policies to approximately maximize the expected reward rate in a loss network. These control policies are easy to implement and perform well both in the transient period as well as in steady state. A major advantage of the penalty approach is that it avoids solving the associated dynamic program. However, a disadvantage of this approach is that it requires the capacity requested by individual requests to be sufficiently small compared to total available capacity. We first solve a related deterministic linear program (LP) and then translate an optimal solution of the LP into an admission control policy for the loss network via an exponential penalty function. We show that the penalty policy is a target-tracking policy—it performs well because the optimal solution of the LP is a good target. We demonstrate that the penalty approach can be extended to track arbitrarily defined target sets. Results from preliminary simulation studies are included.

**1. Introduction.** We consider the following dynamic stochastic allocation problem (details in Section 2). The stochastic system consists of a network of resources (facilities), each with a known fixed capacity. Requests for using this network belong to a diverse set of request classes, differing in the arrival rate, the service duration, the resource requirements and the willingness to pay. There is no waiting room (queue), therefore an arriving request must be either admitted into the system for service and assigned an appropriate resource allocation or rejected (lost) at the instant it arrives. An admitted request occupies the allocated resources for the service duration and releases all the resources simultaneously. The objective of the

Received February 2002; revised September 2003.
[1] Supported in part by NSF Grants CCR-00-09972 and DMS-01-04282.
[2] Supported in part by NSF Grant DMI-01-15034.

system controller is to design an admission control policy that optimizes an appropriate performance measure of the revenue generated.

The stochastic model detailed above is known as a *loss network*. Loss networks model a wide variety of applications where a diverse user population shares a limited collection of resources, for example, telephone networks, local area networks, multiprocessor interconnection architectures, data base structures, mobile radio and broadband packet networks [see Ott and Krishnan (1992), Hui (1990), Kelly (1985), Lagarias, Odlyzko and Zagier (1985), Mitra and Weinberger (1987) and Mitra, Morrison and Ramakrishnan (1996), for details]. Kelly (1991) gave an excellent review of results for loss networks. For a discussion of a related model with loss queues in series, see Ku and Jordan (1997).

A loss network with a single resource is known as a *stochastic knapsack* [Ross and Tsang (1989b)]. Optimality results have been obtained for several restricted classes of admissible policies: complete partitioning policies [Ross and Tsang (1989b)], coordinate convex policies [Foschini and Gopinath (1983), Ross and Tsang (1989b) and Jordan and Varaiya (1994)] and restricted complete sharing policy [Gavois and Rosberg (1994)]. Ross and Yao (1990) discussed monotonicity properties for the stochastic knapsack. See Ross (1995) for a summary of these results.

When capacity requests and service durations of all the request classes are identical, the optimal policy for the stochastic knapsack problem has the following simple form: Accept class $i$ requests if there are at least $\delta_i$ units of capacity free. Such a policy is called a *trunk reservation* policy and the parameters $\delta_i$ are called trunk reservation parameters. This result was established by Miller (1969) [see also Lippman and Ross (1971)]. Several approaches to compute (approximately) optimal trunk reservation parameters $\delta_i$ were discussed by Key (1990), Bean, Gibbens and Zachary (1995) and Reiman and Schwartz (2001). Trunk reservation policies are not optimal when the capacity request or service duration is class dependent [Ross and Tsang (1989a)] nor are they optimal for networks [Key (1990)]. The asymptotic optimality of trunk reservation policies under a limiting regime where the arrival rates and capacity increase together, the Halfin–Whitt regime [Halfin and Whitt (1981)], was established by Hunt and Laws (1993, 1997). For asymptotic optimality results under different limiting regimes, see Kelly (1991), Hunt and Kurtz (1994) and Key (1994).

The optimal capacity allocation problem has also been extensively studied in the revenue management literature. For a recent overview, see McGill and van Ryzin (1999). Unlike the model introduced here, capacity allocation models in the revenue management literature typically assume that there is a finite time horizon over which the capacity must be allocated and that capacity once allocated never becomes available again. Our model is closer to that developed by Savin, Cohen, Gans and Katalan (2000) in the context of the rental industry.

In all previous works on related stochastic allocation models, the associated optimization problem is formulated as a dynamic program (DP), and the optimal policy is the solution of the associated Bellman equation. However, solving the Bellman equation quickly becomes computationally intractable and is, in many cases, EXP-complete [Papadimitriou and Tsitsiklis (1999) and Blondel and Tsitsiklis (2000)]. In practice, therefore, the DP formulation is only used to characterize certain qualitative structural properties of the optimal policy, which then form the basis for heuristic approaches for solving the problem. Optimal DP policies are very sensitive to the time horizon of the problem. Due to end-effects, the optimal DP policies that correspond to different time horizons are usually not compatible. Also, there is no guarantee that steady state optimal policies [e.g., the independent thinning policy; Kelly (1991)], will perform well in the transient period.

In this article, we explore alternative simpler techniques for characterizing approximately optimal policies. We replace the stochastic optimization problem by a suitably constructed linear program (LP). The optimal solution of this LP yields a target point that is translated into an admission control policy using an exponential penalty function. We show that this policy is approximately optimal in the limit where individual resource requests are small compared to the total capacity [Halfin and Whitt (1981)]. Moreover, we show that this penalty policy performs well in the transient period as well.

Our penalty-based approach builds on several disparate research ideas: convex programming bounds for stochastic problems [Gibbens and Kelly (1995), Bertsimas, Paschalidis and Tsitsiklis (1994), Bertsimas and Niño Mora (1999a, b) and Bertsimas and Chryssikou (1999)], asymptotically optimal policies for control and scheduling problems via "fluid" relaxations [Maglaras (2000), Bertsimas and Sethuraman (2002) and Bertsimas, Sethuraman and Gamarnik (2003)] and exponential penalty-based approximation algorithms for linear programming [Shahronki and Matula (1990), Plotkin, Shmoys and Tardos (1991) and Bienstock (2002)]. Exponential penalty functions have also proved useful for admission control and load balancing in an adversarial setting [Aspnes, Azar, Plotkin and Waarts (1997), Azar, Kalyanasundaram, Plotkin, Pruhs and Waarts (1997) and Kamath, Palmon and Plotkin (1998)]. Of this, Kamath, Palmon and Plotkin (1998) is the most relevant to the discussion here.

The summary of our contributions in this article is as follows:

(i) We develop explicit upper bounds for the maximum achievable revenue rate for any time $t \geq 0$. This extends the analysis in Gibbens and Kelly (1995).

(ii) We construct an exponential penalty-based admission control policy that is provably approximately optimal for all times $t \geq 0$ in the Halfin–Whitt limiting regime [Halfin and Whitt (1981)]. The policy is a simple threshold-type policy in an expanded state space. Preliminary simulation studies (see Section 3.4) suggest that the state space expansion is the key to the success of the penalty policy.

(iii) We demonstrate that our approach can be extended to track arbitrary polyhedral target sets.

The organization of this article is as follows. In Section 2 we formulate the admission control problem for a loss network. The framework is Markovian, that is, the arrivals are Poisson and service times are exponentially distributed. In Section 3 we study the single resource model and its various variants. Section 3.4 contains simulation results for this special case and Section 3.5 extends some of the results to the case of general service time distributions. In Section 4 we extend the single-resource results to the network problem. Section 5 presents an extension to control problems where the objective is to ensure that the state of the network lies in a specified target set. Section 6 has some concluding comments and discussion.

**2. Admission control in loss networks.** The stochastic system under consideration consists of a network of $s$ resources (facilities) with capacity $\mathbf{b} \in \mathbf{R}_+^s$, where $b(k) \geq 0$ is the capacity of resource $k = 1, \ldots, s$. Requests for using this network belong to $m$ independent Poisson arrival classes. Class $i$ requests have an arrival rate $\lambda_i$ and a service duration $S_i \sim \exp(\mu_i)$; that is, $S_i$ is exponentially distributed with rate $\mu_i$ (with the exception of Section 3.5). Class $i$ requests are willing to accept any capacity allocation from the set $\mathcal{B}_i = \{\mathbf{b}_{i1}, \ldots, \mathbf{b}_{il_i}\}$, $\mathbf{b}_{ij} \in \mathbf{R}_+^s$, and pay $r_i$ per unit time for the (random) service duration $S_i$. There is no waiting room in the system; therefore, each arriving class $i$ request must either be accepted and admitted into the system (i.e., assigned an admissible capacity allocation $\mathbf{b}_{ij} \in \mathcal{B}_i$) or be rejected at the instant it arrives. When an accepted request departs after service completion, it releases all the allocated resources simultaneously.

We assume that the system is initially empty, that is, $\mathbf{x}(0^-) = \mathbf{0}$ (see Remark 1 in Section 3.1 for a discussion on nonzero initial states). Let $x_{ij}(t)$ denote the number of class $i$ requests currently in the system that are assigned to the allocation $\mathbf{b}_{ij} \in \mathcal{B}_i$. Define $\mathbf{x}_i(t) = (x_{i1}(t), \ldots, x_{il_i}(t)) \in \mathbf{Z}_+^{l_i}$ and $\mathbf{x}(t) = (\mathbf{x}_1(t), \ldots, \mathbf{x}_m(t)) \in \mathbf{Z}_+^l$, where $l = \sum_{i=1}^m l_i$. A request of class $i$ can be assigned a capacity allocation $\mathbf{b}_{ij}$ only if there is sufficient capacity to accommodate it, that is,

$$(1) \qquad \sum_{i'=1}^m \sum_{j'=1}^{l_{i'}} x_{i'j'}(t)\mathbf{b}_{i'j'} + \mathbf{b}_{ij} \leq \mathbf{b},$$

where the inequality is interpreted component by component. The system controller is permitted to reject requests even if there is sufficient capacity to accommodate them. The instantaneous reward rate $R(t)$ at time $t$ is given by

$$(2) \qquad R(t) = \sum_{i=1}^{m} r_i\left(\sum_{j=1}^{l_i} x_{ij}(t)\right) = \sum_{i=1}^{m} r_i(\mathbf{1}^T \mathbf{x}_i(t)).$$

This stochastic model is called a *loss network* [Kelly (1991)].

Let $T_{(i,n)}$, $i = 1, \ldots, m$, $n \geq 1$, denote the arrival epoch of the $n$th class $i$ request. Since all admission decisions are made at arrival epochs, a feasible admission control policy $\boldsymbol{\pi}$ is described as follows:

(a) A policy $\boldsymbol{\pi}$ is a collection of random variables $\boldsymbol{\pi} = \{\pi_{(i,n)} : i = 1, \ldots, m, n \geq 1\}$, with $\pi_{(i,n)} \in \{0, 1, \ldots, l_i\}$, where $\pi_{(i,n)} = 0$ denotes that class $i$ request arriving at the epoch $T_{(i,n)}$ is rejected and $\pi_{(i,n)} = j$ ($\geq 1$) denotes that the request is assigned to $\mathbf{b}_{ij} \in \mathcal{B}_i$.

(b) The random variable $\pi_{(i,n)}$ is measurable with respect to the $\sigma$-algebra generated by the past arrival epochs $\{T_{(p,q)} : p = 1, \ldots, m, q \geq 1, T_{(p,q)} \leq T_{(i,n)}\}$, the past actions $\{\pi_{(p,q)} : p = 1, \ldots, m, q \geq 1, T_{(p,q)} \leq T_{(i,n)}\}$ and the state process $\{\mathbf{x}^\pi(t) : t \leq T_{(i,n)}\}$, where the notation $\mathbf{x}^\pi$ emphasizes that the state process is itself a function of past actions.

(c) The state process $\{\mathbf{x}^\pi(t) : t \geq 0\}$ does not violate capacity constraints, that is, $\sum_{i=1}^{m} \sum_{j=1}^{l_i} x_{ij}^\pi(t) \mathbf{b}_{ij} \leq \mathbf{b}$ for all $t \geq 0$. (Rejection is the only feasible action when adequate capacity is not available.)

Let $R^\pi(t) = \sum_{i=1}^{m} r_i(\mathbf{1}^T \mathbf{x}_i^\pi(t))$ denote the instantaneous reward rate of the policy $\boldsymbol{\pi}$ at time $t$. The objective of the controller is to choose a feasible policy $\boldsymbol{\pi}$ that maximizes some performance measure on the reward rate process $\{R^\pi(t) : t \geq 0\}$. Appropriate performance measures for finite time horizon problems are either expected total reward $\mathbf{E}[\int_0^T R^\pi(s) \, ds]$ or expected discounted reward $\mathbf{E}[\int_0^T e^{-\beta s} R^\pi(s) \, ds]$, $\beta > 0$; for the infinite time horizon problems, the appropriate measures are either expected discounted reward $\mathbf{E}[\int_0^\infty e^{-\beta s} R^\pi(s) \, ds]$, $\beta > 0$, or long-run average reward $\lim_{T \to \infty} \frac{1}{T} \mathbf{E}[\int_0^T R^\pi(s) \, ds]$.

As mentioned in Section 1, our goal is to construct feasible policies that perform well both in the transient period as well as in steady state. We first establish an upper bound $R^*(t)$ on the achievable expected reward rate $\mathbf{E}[R^\pi(t)]$ and then construct a feasible policy $\bar{\boldsymbol{\pi}}$ with expected reward rate $\mathbf{E}[\bar{R}(t)] \approx R^*(t)$. Thus, the policy $\bar{\boldsymbol{\pi}}$ satisfies

$$\mathbf{E}\left[\int_0^T e^{-\beta s} R^\pi(s) \, ds\right] \leq \int_0^T e^{-\beta s} R^*(s) \, ds \approx \mathbf{E}\left[\int_0^T e^{-\beta s} \bar{R}(s) \, ds\right], \qquad \beta \geq 0,$$

that is, the policy $\bar{\boldsymbol{\pi}}$ is approximately optimal for any finite time horizon, and

$$\lim_{t\to\infty}\frac{1}{T}\mathbf{E}\bigg[\int_0^T R^\pi(s)\,ds\bigg] \leq \lim_{T\to\infty}\frac{1}{T}\int_0^T R^*(s)\,ds \approx \lim_{T\to\infty}\frac{1}{T}\mathbf{E}\bigg[\int_0^T \bar{R}(s)\,ds\bigg],$$

that is, the policy $\bar{\boldsymbol{\pi}}$ is approximately optimal in the steady state as well.

**3. Single-resource model.** This section focuses on the loss network with $s=1$ (i.e., the stochastic knapsack). The details of the single-resource model are as follows. The system is assumed to be initially empty [i.e. $\mathbf{x}(0^-)=\mathbf{0}$]. Requests belong to $m$ Poisson arrival classes. Request class $i$ has arrival rate $\lambda_i$, capacity request $b_i$ (without loss of generality, one can assume that the set $\mathcal{B}_i$ is a singleton), service duration $S_i \sim \exp(\mu_i)$, and reward rate $r_i$ per unit time. All the requests arrive at a common resource with capacity $b \in (0,\infty)$. There is no waiting space (queue); therefore, each arriving request must either be admitted into service or rejected at the instant it arrives [see Cosyn and Sigman (2004) and Cosyn (2003) for extensions to queues]. Requests may be rejected even if there was adequate capacity available.

Note that if the total capacity $b$ is an integer and $b_i = 1$, $1 \leq i \leq m$, then $b$ can be identified as the number of servers in a standard queuing model. In particular, if requests are always served when capacity exists, then this is simply an $M/M/b$ loss queue. Thus, it helps to imagine that each accepted request has its own server. In this light, the loss network introduced in Section 2 can be viewed as a collection of such server models, all working together in parallel.

The layout of this section is as follows. In Section 3.1 we develop an upper bound on the achievable reward rate. In Section 3.2 we construct an approximately optimal penalty-based policy. Section 3.3 investigates the penalty policy in the Halfin–Whitt limiting regime [Halfin and Whitt (1981)]. In Section 3.4 we simulate the transient behavior of the proposed control policy and compare its performance to thinning policies introduced by Kelly (1991). Section 3.5 discusses the extension to general service times.

3.1. *Upper bound on the achievable reward rate.* Let $\boldsymbol{\pi}$ denote any feasible control policy for the single-resource model. Let $x_i^\pi(t)$ denote the number of the class $i$ requests in service at time $t$. Since feasibility implies that $\sum_{i=1}^m b_i x_i^\pi(t) \leq b$, we have

$$\sum_{i=1}^m b_i \mathbf{E}[x_i^\pi(t)] \leq b. \tag{3}$$

Moreover, $\mathbf{E}[x_i^\pi(t)] \leq \mathbf{E}[q_i(t)]$, where $q_i(t)$ is the number of class $i$ requests as time $t$ in an infinite capacity system with no admission control. Recall that we assume that the system is initially empty, therefore [see, e.g., page 75 in Wolff (1989)], $\mathbf{E}[q_i(t)] = \rho_i(1 - \exp(-\mu_i t))$. Hence,

$$\boldsymbol{\alpha} = \left(\frac{1}{\rho_1}\mathbf{E}[x_1^\pi(t)], \ldots, \frac{1}{\rho_m}\mathbf{E}[x_m^\pi(t)]\right)$$

is feasible for the linear program

(4)
$$\begin{aligned}\text{maximize} \quad & \sum_{i=1}^n r_i \rho_i \alpha_i \\ \text{subject to} \quad & \sum_{i=1}^m b_i \rho_i \alpha_i \leq b, \\ & 0 \leq \alpha_i \leq 1 - \exp(-\mu_i t), \qquad i = 1, \ldots, m.\end{aligned}$$

Let $\boldsymbol{\alpha}^*(t)$ denote an optimal solution and let $R^*(t)$ denote the optimal value of (4). Then

(5) $$\mathbf{E}[R^\pi(t)] = \sum_{i=1}^m r_i \rho_i \left(\frac{1}{\rho_i}\mathbf{E}[x_i^\pi(t)]\right) \leq R^*(t).$$

In the next section we propose a policy that controls the system by penalizing deviations from a desired target state. From (4) and (5), it follows that for a policy $\boldsymbol{\pi}$ to be approximately optimal, the expected number $\mathbf{E}[x_i^\pi(t)]$ of accepted class $i$ requests must be approximately $x_i^*(t) = \alpha_i^*(t)\rho_i$. Thus, $\mathbf{x}^*(t) = (x_1^*(t), \ldots, x_m^*(t))^T$ would be the natural target state for the penalty policy. Unfortunately we are only able to establish that a penalty policy can successfully track a fixed target. The natural fixed target is $x_i^* = \alpha_i^* \rho_i$, $i = 1, \ldots, m$, where $\boldsymbol{\alpha}^* = (\alpha_1, \ldots, \alpha_m)^T$ is an optimal solution of the "steady state" analog of (4):

(6)
$$\begin{aligned}\text{maximize} \quad & \sum_{i=1}^n r_i \rho_i \alpha_i \\ \text{subject to} \quad & \sum_{i=1}^m b_i \rho_i \alpha_i \leq b, \\ & 0 \leq \alpha_i \leq 1, \qquad i = 1, \ldots, m.\end{aligned}$$

Let $R^*$ denote the optimal value of (6). Next, we bound $R^*(t)$ in terms of the steady state quantities $\boldsymbol{\alpha}^*$, $R^*$ and the problem parameters. Since $\boldsymbol{\alpha}$ feasible for (4) must satisfy $\alpha_i \leq 1 - e^{-\mu_i t}$, $i = 1, \ldots, m$, it follows that

(7) $$R^*(t) \leq \sum_{i=1}^m r_i \rho_i (1 - \exp(-\mu_i t)).$$

The linear programming dual of (4) is

$$
\begin{aligned}
\text{minimize} \quad & ub + \sum_{i=1}^{m} v_i(1 - \exp(-\mu_i t)) \\
\text{subject to} \quad & v_i + b_i \rho_i u \geq r_i \rho_i, \quad i = 1, \ldots, m, \\
& \mathbf{v} \geq \mathbf{0}, \quad u \geq 0.
\end{aligned}
\tag{8}
$$

Taking the limit $t \to \infty$ in (8) we get the dual of the steady state LP (6):

$$
\begin{aligned}
\text{minimize} \quad & ub + \mathbf{1}^T \mathbf{v} \\
\text{subject to} \quad & v_i + b_i \rho_i u \geq r_i \rho_i, \quad i = 1, \ldots, m, \\
& \mathbf{v} \geq \mathbf{0}, \quad u \geq 0.
\end{aligned}
\tag{9}
$$

Let $(u^*, \mathbf{v}^*)$ denote any optimal solution of (9), $U = \{i : \alpha_i^* = 1\}$ and $U^c = \{i : i \notin U\}$. Then it follows that

$$
R^*(t) \leq u^* b + \sum_{i=1}^{m} v_i^*(1 - \exp(-\mu_i t)) \tag{10}
$$

$$
= \sum_{i=1}^{m} r_i \rho_i \alpha_i^* - \sum_{i \in U} v_i^* \exp(-\mu_i t) \tag{11}
$$

$$
= \sum_{i=1}^{m} r_i \rho_i \alpha_i^* - \sum_{i \in U} (r_i \rho_i - b_i \rho_i u^*) \alpha_i^* \exp(-\mu_i t) \tag{12}
$$

$$
= \sum_{i=1}^{m} r_i \rho_i \alpha_i^* (1 - \exp(-\mu_i t)) + u^* \left( \sum_{i=1}^{m} b_i \rho_i \alpha_i^* \exp(-\mu_i t) \right) \tag{13}
$$

$$
\leq \sum_{i=1}^{m} r_i \rho_i \alpha_i^* (1 - \exp(-\mu_i t)) + u^* b \exp(-\mu_{\min} t), \tag{14}
$$

where (10) is implied by the fact that $(u^*, \mathbf{v}^*)$ is feasible for the dual LP (8); (11)–(13) all follow from complementary slackness conditions [Luenberger (1984)]; and $\mu_{\min} = \min_{1 \leq i \leq m}\{\mu_i\}$. From (7) and (14) we have the following result.

THEOREM 1.  *The reward rate $R^\pi(t)$ of any feasible policy $\pi$ satisfies*

$$
\begin{aligned}
\mathbf{E}[R^\pi(t)] &\leq R^*(t) \\
&\leq \min\Bigg\{ \sum_{i=1}^{m} r_i \rho_i (1 - \exp(-\mu_i t)), \\
&\qquad\qquad \sum_{i=1}^{m} r_i \rho_i \alpha_i^* (1 - \exp(-\mu_i t)) + u^* b \exp(-\mu_{\min} t) \Bigg\},
\end{aligned}
\tag{15}
$$

where $R^*(t)$ is the optimal value of the LP (4), $\boldsymbol{\alpha}^*$ is an optimal solution of the steady state LP (6) and $(u^*, \mathbf{v}^*)$ is an optimal solution of the steady state dual LP (9).

The first term in the upper bound on $R^*(t)$ is active for $t \leq 1/\mu_{\max}$, where $\mu_{\max} = \max_{1 \leq i \leq m}\{\mu_i\}$, whereas the second is active for $t \geq 1/\mu_{\min}$.

REMARK 1. Although we assume that the system is initially empty, all the results in this article extend to the case where the initial state $\mathbf{x}(0^-) \neq \mathbf{0}$. For example, when $\mathbf{x}(0^-) = \mathbf{x}^0 \neq \mathbf{0}$, the bound analogous to (15) is given by

$$R^*(t) \leq \min\bigg\{\sum_{i=1}^m r_i\rho_i(1 - \exp(-\mu_i t)) + \sum_{i=1}^m r_i x_i^0 \exp(-\mu_i t),$$
$$\sum_{i=1}^m r_i\rho_i\alpha_i^*(1 - \exp(-\mu_i t)) + u^* b \exp(-\mu_{\min} t)$$
$$+ \sum_{i=1}^m \frac{v_i^* x_i^0}{\rho_i} \exp(-\mu_i t)\bigg\}.$$

The results in this section bear close resemblance to the notion of fluid operating points introduced by Harrison (2003). However, unlike the development here, Harrison employed the fluid model only to define a nominal operating point—the control policy is designed using a heavy-traffic limit associated with this operating point.

3.2. *Exponential penalty function and penalty control policy.* Kelly (1991) established that, under fairly general conditions, an independent thinning policy that accepts each incoming class $i$ request with probability $\alpha_i^*$, provided there is enough capacity, approximately optimizes the expected reward rate in steady state. However, for small $t$, thinning underutilizes the capacity and, therefore, the expected reward rate of the thinning policy is significantly smaller than the upper bound (7). Moreover, since thinning only changes the effective arrival rate, it is not able to effectively control the variance of the reward rate. Our goal is to construct a policy that does not suffer from these drawbacks. We will first informally motivate the structure of the policy and then establish its properties rigorously.

Consider the following modification to the original system. Suppose each rejected class $i$ request, instead of immediately leaving the system, is assigned to an alternate infinite capacity server where it lives out its service time and then leaves. [In practice, each time a request is rejected the policy will add one request to the alternate server with a service time $S_i \sim \exp(\mu_i)$.]

From the analysis leading to the LP (4), it follows that for the expected reward rate $\mathbf{E}[R(t)]$ to be close to the bound (15), one requires $\mathbf{E}[x_i(t)] \approx x_i^*(t) = \alpha_i(t)\rho_i$, $i = 1, \ldots, m$. Let $y_i(t)$ denote the number of class $i$ requests in the alternate server at time $t$. Then $\mathbf{E}[x_i(t)] + \mathbf{E}[y_i(t)] = \mathbf{E}[q_i(t)] = \rho_i(1 - \exp(-\mu_i t))$. Thus, an equivalent condition for optimality is that $\mathbf{E}[y_i(t)] \approx y_i^*(t) = \rho_i(1 - \exp(-\mu_i t) - \alpha_i(t))$. Let $\Psi_i(x_i, y_i)$ be a penalty function that penalizes deviations from the desired target state $(x_i^*(t), y_i^*(t))$. Since keeping $(x_i, y_i) \approx (x_i^*(t), y_i^*(t))$ is equivalent to minimizing the penalty function, a control policy that accepts a request, provided there is adequate capacity and $\Psi_i(x_i + 1, y_i) \leq \Psi_i(x_i, y_i + 1)$, may be close to optimal. Such a policy can be thought of as iteratively solving the nonlinear optimization problem $\min_{x,y} \Psi_i(x, y)$ with the added restriction that it can take a step only when there is an arrival and the step length is restricted to 1. Moreover, periodically the state $(x_i, y_i)$ gets perturbed in a uncontrollable manner by requests leaving the system. From related results in the nonlinear optimization literature [see, e.g., Luenberger (1984)], it follows that such a penalty-based control policy is likely to be successful provided the gradient of the penalty $\Psi_i$ is sufficiently "large" around the target state $(x_i^*, y_i^*)$, the step length of 1 is a "small" step in an appropriately defined norm and the frequency of correcting steps is sufficiently higher than the frequency of the perturbing steps (i.e., $\rho_i = \lambda_i/\mu_i \gg 1$). The relationship of penalty function and nonlinear optimization is further discussed in Section 6.

In this article, we use a penalty function of the form

$$\Psi_i(x_i, y_i) = \exp\left(\beta \frac{x_i(t)}{x_i^*(t)}\right) + \exp\left(\beta \frac{y_i(t)}{y_i^*(t)}\right).$$

This choice is motivated by the fact that the exponential function is an eigenfunction of the underlying Markov process and that, for this choice, moment generating functions can be used to characterize the behavior of the penalty policy. Note that although the penalty method can be formulated without any reference to the rejected requests $y_i$, the form that we propose does not permit us to do so. In our penalty function we need $y_i$ to ensure that the number of accepted requests $x_i$ does not drop too low. In the rest of this section, we rigorously establish these informal ideas.

Since we are interested in approximating the upper bound (15), we drop from consideration all those classes with $\alpha_i^* = 0$. As proposed above, we add a fictitious infinite capacity system. We will refer to the original system as system 0 and the fictitious system as system 1. The state of the augmented network at time $t$ is $\mathbf{s}(t) = (\mathbf{x}(t), \mathbf{y}(t)) \in \mathbf{Z}_+^{2m}$. The state vector $\mathbf{x}(t) = (x_1(t), \ldots, x_m(t))$, where $x_i(t)$ is the number of class $i$ requests in system 0 at time $t$, describes the state of system 0. Similarly, $\mathbf{y}(t) = (y_1(t), \ldots, y_m(t))$ describes the state of the fictitious system 1 at time $t$.

The state $\mathbf{s} = (\mathbf{x}, \mathbf{y})$ is assigned a penalty $\Psi(\mathbf{s})$ given by

$$\Psi(\mathbf{s}) = \sum_{i=1} \underbrace{\left[\exp\left(\beta \cdot \frac{b_i x_i}{c_i^0}\right) + \exp\left(\beta \cdot \frac{b_i y_i}{c_i^1}\right)\right]}_{\Psi_i(\mathbf{s}_i)}, \tag{16}$$

where $(\mathbf{c}^0, \mathbf{c}^1) \in \mathbf{R}_+^{2m}$ and $\mathbf{s}_i = (x_i, y_i)$ denotes the components of $\mathbf{s}$ that correspond to class $i$. There are two competing requirements on the multiplier $\beta$: we need $\beta$ to be large to ensure that the penalty function $\Psi(\mathbf{s})$ is sufficiently steep; on the other hand, we also have to ensure that the impact of a single arrival or departure on the penalty value is sufficiently small. The precise bound on $\beta$ is given by (22). The capacities $(\mathbf{c}^0, \mathbf{c}^1)$ determine the "steady-state" target state of the penalty policy. As mentioned previously, we choose a fixed target because we are unable to establish that penalty policies can track time-varying targets. The transient performance is controlled by suitably initializing the fictitious system 1.

The penalty policy $\bar{\boldsymbol{\pi}}$ is defined as follows. Let $\{\bar{\mathbf{s}}(t) = (\bar{\mathbf{x}}(t), \bar{\mathbf{y}}(t)) : t \geq 0\}$ denote the state process under the control $\bar{\boldsymbol{\pi}}$. At time $t = 0^-$, the state of the original system $\mathbf{x}(0^-) = \mathbf{0}$, the state of the fictitious infinite capacity system 1 is initialized to $\bar{\mathbf{y}}(0^-)$ [the precise value of $\bar{\mathbf{y}}(0^-)$ is specified later] and a service time $S_i \sim \exp(\mu_i)$ is generated for each of the $\bar{y}_i(0^-)$ class $i$ requests in system 1, $i = 1, \ldots, m$.

At time $t \geq 0$, an arriving class $i$ request is accepted by the control policy $\bar{\boldsymbol{\pi}}$ (i.e., routed to system 0) provided

$$\frac{\partial \Psi_i(\bar{\mathbf{s}}_i(t))}{\partial x_i} \leq \frac{\partial \Psi_i(\bar{\mathbf{s}}_i(t))}{\partial y_i} \tag{17}$$

and the capacity constraint on system 0 is not violated, that is,

$$\sum_{i'=1}^{m} b_{i'} \bar{x}_{i'}(t) + b_i \leq b; \tag{18}$$

otherwise it is rejected (i.e., routed to system 1) and the policy $\bar{\boldsymbol{\pi}}$ attaches to it a service time $S_i \sim \exp(\mu_i)$ independent of everything else. Since the admission condition (17) is equivalent to

$$\frac{\bar{x}_i(t)}{c_i^0} \leq \frac{\bar{y}_i(t)}{c_i^1} + \frac{1}{\beta b_i} \log\left(\frac{c_i^0}{c_i^1}\right), \tag{19}$$

it is clear that the policy $\bar{\boldsymbol{\pi}}$ is a threshold-type policy in the expanded state space $\mathbf{s} = (\mathbf{x}, \mathbf{y}) \in \mathbf{Z}_+^{2m}$.

The capacities $(\mathbf{c}^0, \mathbf{c}^1)$, the parameter $\beta$ and the initial state $\bar{\mathbf{y}}(0^-)$ are defined in terms of a perturbation parameter $\varepsilon \in (0, \frac{1}{4})$. Define an $\varepsilon$-perturbation of the steady state LP (6) as

(20)
$$\begin{aligned}
\text{maximize} \quad & \sum_{i=1}^{m} r_i \rho_i \alpha_i \\
\text{subject to} \quad & \sum_{i=1}^{m} b_i \rho_i \alpha_i \leq \frac{b}{1+4\varepsilon}, \\
& 0 \leq \alpha_i \leq 1, \qquad i=1,\ldots,m.
\end{aligned}$$

Let $\boldsymbol{\alpha}^\varepsilon$ denote an optimal solution of this perturbed LP (20). Then the capacities $(\mathbf{c}^0, \mathbf{c}^1)$ are given by

(21) $\quad c_i^0 = (1+4\varepsilon)\alpha_i^\varepsilon b_i \rho_i, \qquad c_i^1 = (1+4\varepsilon)(1-\alpha_i^\varepsilon)b_i \rho_i, \qquad i=1,\ldots,m,$

and $\beta$ must satisfy

(22) $\qquad \beta \leq \varepsilon \min\left\{ \min_{\{i:\, 1\leq i \leq m\}}\left\{\frac{c_i^0}{b_i}\right\}, \min_{\{i:\, i \in U_\varepsilon^c\}}\left\{\frac{c_i^1}{b_i}\right\} \right\}$

(23) $\qquad = \varepsilon(1+4\varepsilon) \min\left\{ \min_{\{i:\, 1\leq i \leq m\}}\{\alpha_i^\varepsilon \rho_i\}, \min_{\{i:\, i \in U_\varepsilon^c\}}\{(1-\alpha_i^\varepsilon)\rho_i\} \right\},$

where $U_\varepsilon^c = \{i : \alpha_i^\varepsilon < 1,\, i=1,\ldots,m\}$. The bound (22) formalizes the notion that the change in the penalty value associated with a single arrival or departure must be small [the bounds (22) and (23) are identical]. Since parameter $\beta$ must be sufficiently large for the penalty policy to perform well, the bound (23) implies that penalty policy is likely to perform well when the incoming load $\rho_i \gg 1$. Although the request sizes $b_i$ are not explicitly present, the bounds (22) and (23) impose an implicit upper bound on the $b_i$'s via the capacity constraint $\sum_i b_i \rho_i \alpha_i \leq b$.

We establish a lower bound on the expected reward rate $\mathbf{E}[\bar{R}(t)]$ of the policy $\bar{\boldsymbol{\pi}}$ by comparing it to a related infeasible policy $\tilde{\boldsymbol{\pi}}$. The policy $\tilde{\boldsymbol{\pi}}$ is identical to $\bar{\boldsymbol{\pi}}$ except that it does not respect the system 0 capacity constraints; that is, the policy $\tilde{\boldsymbol{\pi}}$ routes an incoming class $i$ request to system 0 whenever

(24) $$\frac{\partial \Psi_i(\tilde{\mathbf{s}}_i(t))}{\partial x_i} \leq \frac{\partial \Psi_i(\tilde{\mathbf{s}}_i(t))}{\partial y_i},$$

where $\{\tilde{\mathbf{s}}(t) = (\tilde{x}(t), \tilde{y}(t)) : t \geq 0\}$ denotes the state process that corresponds to the policy $\tilde{\boldsymbol{\pi}}$. Since the various request classes interact only through the capacity constraints, the policy $\tilde{\boldsymbol{\pi}}$ controls each class independently.

We establish a bound on the total derivative $(d/dt)\mathbf{E}[\Psi(\tilde{\mathbf{s}}(t))]$, which implies that if the initial state $\tilde{\mathbf{y}}(0^-)$ is suitably chosen, the penalty $\mathbf{E}[\Psi(\tilde{\mathbf{s}}(t))]$ is a uniformly bounded function of time.

LEMMA 1. *Suppose $\varepsilon < \frac{1}{4}$, $(\mathbf{c}^0, \mathbf{c}^1)$ are given by (21) and $\beta$ satisfies (22). Then, for all $i = 1, \ldots, m$ and $t \geq 0$,*

$$\frac{d}{dt}\mathbf{E}[\Psi_i(\tilde{\mathbf{s}}_i(t))] \leq (1-\varepsilon)\mu_i(2e^{(1-\varepsilon/2)\beta} - \mathbf{E}[\Psi_i(\tilde{\mathbf{s}}_i(t))]).$$

PROOF. Fix a request class $i$. Define $\mathbf{E}_t[\Psi_i(\tilde{\mathbf{s}}_i(u))] = \mathbf{E}[\Psi_i(\tilde{\mathbf{s}}_i(u)) \mid \mathcal{F}_t]$, $u \geq t$, where $\mathcal{F}_t$ is the filtration generated by events up to $t$. Then

$$\frac{d}{dt}\mathbf{E}_t[\Psi_i(\tilde{\mathbf{s}}(t))] = \mathcal{A}\Psi_i(\tilde{\mathbf{s}}(t)),$$

where $\mathcal{A}$ is the generator of the stochastic process $\{\tilde{\mathbf{s}}(t) : t \geq 0\}$. Let $\tilde{\pi}_i(t)$ denote the routing decision of the policy $\tilde{\boldsymbol{\pi}}$ at time $t$, that is,

$$\tilde{\pi}_i(t) = \begin{cases} 1, & \frac{\partial \Psi_i}{\partial x_i} \leq \frac{\partial \Psi_i}{\partial y_i}, \\ 0, & \text{otherwise.} \end{cases}$$

Then

$$\mathcal{A}\Psi_i(\tilde{\mathbf{s}}(t)) = \lambda_i[(\Psi_i(\tilde{x}_i + \tilde{\pi}_i(t), \tilde{y}_i) - \Psi_i(\tilde{x}_i, \tilde{y}_i))$$
$$+ (\Psi_i(\tilde{x}_i, \tilde{y}_i + (1 - \tilde{\pi}_i(t))) - \Psi_i(\tilde{x}_i, \tilde{y}_i))]$$
$$+ \mu_i[x_i(\Psi_i(\tilde{x}_i - 1, \tilde{y}_i) - \Psi_i(\tilde{x}_i, \tilde{y}_i))$$
$$+ y_i(\Psi_i(\tilde{x}_i, \tilde{y}_i - 1) - \Psi_i(\tilde{x}_i, \tilde{y}_i))],$$

where we have suppressed the time dependence of $(\tilde{x}_i, \tilde{y}_i)$. From the Taylor series expansion, it follows that $e^x \leq 1 + x + x^2$ for all $|x| \leq 1$ and from the bound (22) we have that $\max\{\beta b_i/c_i^0, \beta b_i/c_i^1\} \leq \varepsilon$. Therefore,

$$\mathcal{A}\Psi_i(\tilde{\mathbf{s}}(t)) \leq (1+\varepsilon)\mu_i\left(\frac{\partial \Psi_i}{\partial x_i} \cdot \tilde{\pi}_i(t)\rho_i + \frac{\partial \Psi_i}{\partial y_i} \cdot (1 - \tilde{\pi}_i(t))\rho_i\right)$$
$$- (1-\varepsilon)\mu_i\left(\frac{\partial \Psi_i}{\partial x_i} \cdot \tilde{x}_i(t) + \frac{\partial \Psi_i}{\partial y_i} \cdot \tilde{y}_i(t)\right).$$

Since $\tilde{\pi}_i(t)$ minimizes the increase in penalty, it follows that

$$\frac{\partial \Psi_i}{\partial x_i} \cdot \tilde{\pi}_i(t)\rho_i + \frac{\partial \Psi_i}{\partial y_i} \cdot (1 - \tilde{\pi}_i(t))\rho_i \leq \frac{\partial \Psi_i}{\partial x_i} \cdot x_i^\varepsilon + \frac{\partial \Psi_i}{\partial y_i} \cdot y_i^\varepsilon$$

for any $x_i^\varepsilon + y_i^\varepsilon = \rho_i$, $x_i^\varepsilon, y_i^\varepsilon \geq 0$. In particular, choose

(25) $$x_i^\varepsilon = \alpha_i^\varepsilon \rho_i, \qquad y_i^\varepsilon = (1 - \alpha_i^\varepsilon)\rho_i.$$

Then, we have

$$\mathcal{A}\Psi_i(\tilde{\mathbf{s}}(t)) \leq (1+\varepsilon)\mu_i\left(\frac{\partial \Psi_i}{\partial x_i} \cdot x_i^\varepsilon + \frac{\partial \Psi_i}{\partial y_i} \cdot y_i^\varepsilon\right)$$

$$- (1-\varepsilon)\mu_i\left(\frac{\partial \Psi_i}{\partial x_i}\cdot \tilde{x}_i(t) + \frac{\partial \Psi_i}{\partial y_i}\cdot \tilde{y}_i(t)\right)$$

$$= (1-\varepsilon)\mu_i\left[\frac{\partial \Psi_i}{\partial x_i}\left(\frac{(1+\varepsilon)}{(1-\varepsilon)}x_i^\varepsilon - \tilde{x}_i(t)\right) + \frac{\partial \Psi_i}{\partial y_i}\left(\frac{(1+\varepsilon)}{(1-\varepsilon)}y_i^\varepsilon - \tilde{y}_i(t)\right)\right]$$

(26) $$\leq (1-\varepsilon)\mu_i\left[\Psi_i\left(\frac{(1+\varepsilon)}{(1-\varepsilon)}\mathbf{s}_i^\varepsilon\right) - \Psi_i(\tilde{\mathbf{s}})\right]$$

(27) $$\leq (1-\varepsilon)\mu_i[\Psi_i((1+3\varepsilon)\mathbf{s}_i^\varepsilon) - \Psi_i(\tilde{\mathbf{s}})],$$

where (26) follows from the convexity of $\Psi_i$ and (27) holds because $\frac{1+\varepsilon}{1-\varepsilon} \leq 1 + 3\varepsilon$ for all $\varepsilon < \frac{1}{4}$. From (21) and (25), it follows that $(1+3\varepsilon)\max\{b_i x_i^\varepsilon/c_i^0, b_i y_i^\varepsilon/c_i\} = \frac{1+3\varepsilon}{1+4\varepsilon} \leq 1 - \frac{\varepsilon}{2}$. Consequently,

$$\frac{d}{dt}\mathbf{E}_t[\Psi_i(\tilde{\mathbf{s}}_i(t))] \leq (1-\varepsilon)\mu_i[2e^{(1-\varepsilon/2)\beta} - \Psi_i(\tilde{\mathbf{s}}(t))].$$

The result can now be concluded from the Lebesgue bounded convergence theorem by recognizing that for all sufficiently close $s \geq t$, $(\mathbf{E}_t[\Psi_i(\tilde{\mathbf{s}}_i(s))] - \Psi_i(\tilde{\mathbf{s}}_i(t)))/(s-t)$ can be bounded above by a fixed random variable. □

LEMMA 2. *Suppose $\varepsilon < \frac{1}{4}$, $(\mathbf{c}^0, \mathbf{c}^1)$ are given by (21), $\beta$ satisfies (22) and the initial state $\tilde{\mathbf{s}}_i(0^-) = (\mathbf{0}, \tilde{\mathbf{y}}(0^-))$ satisfies $\Psi_i(\tilde{\mathbf{s}}_i(0^-)) \leq 2\exp((1-\varepsilon/2)\beta)$, $i = 1,\ldots,m$. Then, for all $i = 1,\ldots,m$ and $t \geq 0$,*

(28) $$\mathbf{E}[\Psi_i(\tilde{\mathbf{s}}_i(t))] \leq 2e^{(1-\varepsilon/2)\beta}.$$

PROOF. Fix a request class $i$. Suppose the conclusion does not hold. Define $f_i(t) = \mathbf{E}[\Psi_i(\tilde{\mathbf{s}}_i(t))]$ and $f^* = 2\exp((1-\varepsilon/2)\beta)$. Then Lemma 1 implies that $\frac{df(t)}{dt} \leq (1-\varepsilon)\mu_i(f^* - f_i(t))$.

Let $\tau$ be any time instant when $f(\tau) > f^*$. Since $f(t)$ is a continuous function of $t$ and $f(0^-) \leq f^*$, there exists $s < \tau$ such that $f(s) = f^*$ and $f(t) \geq f^*$ for all $s \leq t \leq \tau$. By construction, $f(\tau) > f^* = f(s)$, but by the fundamental theorem of calculus, we have

$$f(\tau) - f(s) = \int_s^\tau \frac{df(u)}{du}\,du \leq \int_s^\tau (1-\varepsilon)\mu_i(f^* - f(u))\,du \leq 0,$$

a contradiction. □

The bound (28) implies the following results.

LEMMA 3. *Suppose $\varepsilon < \frac{1}{4}$, $(\mathbf{c}^0, \mathbf{c}^1)$ are given by (21) and $\beta$ satisfies (22).*

(i) *Let $\tilde{w}(t) = \sum_{i=1}^m b_i \tilde{x}_i(t)$ and suppose $\Psi_i(\tilde{\mathbf{s}}_i(0^-)) \leq 2\exp((1-\varepsilon/2)\beta)$, $i = 1,\ldots,m$. Then*

(29) $$\mathbf{E}[(\tilde{w}(t) - b)^+] \leq (1+4\varepsilon) \cdot \frac{2e^{-\varepsilon\beta/2}}{\beta} \cdot b.$$

(ii) *Suppose* $\tilde{y}_i(0^-) = (1 - \alpha_i^\varepsilon)\rho_i$, $i = 1, \ldots, m$. *Then the reward rate* $\tilde{R}(t)$ *of the policy* $\tilde{\pi}$ *satisfies*

$$\text{(30)} \qquad \mathbf{E}[\tilde{R}(t)] \geq \sum_{i=1}^{m} \alpha_i^\varepsilon r_i \rho_i (1 - e^{-\mu_i t}) - \zeta \sum_{i=1}^{m} (1 - \alpha_i^\varepsilon) r_i \rho_i,$$

*where* $\boldsymbol{\alpha}^\varepsilon$ *is an optimal solution of the perturbed LP* (20) *and*

$$\zeta = \left(\frac{\log(2)}{\beta} + 1 - \frac{\varepsilon}{2}\right)(1 + 4\varepsilon) - 1.$$

PROOF. Let $V_t = \{\tilde{w}(t) = \sum_{i=1}^{m} b_i \tilde{x}_i(t) > b\}$. Then

$$\text{(31)} \quad \exp\left(\frac{\beta}{b} \cdot \mathbf{E}[(\tilde{w}(t) - b)^+]\right) \leq \mathbf{E}\left[\exp\left(\frac{\beta}{b}(\tilde{w}(t) - b)^+\right)\right]$$

$$= \mathbf{P}(V_t^c) + \mathbf{E}\left[\exp\left(\frac{\beta}{b} \cdot (\tilde{w}(t) - b)\right); V_t\right]$$

$$\leq 1 + \mathbf{E}\left[\exp\left(\frac{\beta}{b} \cdot (\tilde{w}(t) - b)\right)\right]$$

$$\text{(32)} \qquad\qquad = 1 + e^{-\beta} \prod_{1 \leq i \leq m} \mathbf{E}\left[\exp\left(\beta \cdot \frac{b_i \tilde{x}_i(t)}{b}\right)\right],$$

where (31) follows from Jensen's inequality. Moreover,

$$\mathbf{E}\left[\exp\left(\beta \cdot \frac{b_i \tilde{x}_i(t)}{b}\right)\right] = \mathbf{E}\left[\left(\exp\left(\beta \cdot \frac{b_i x_i(t)}{c_i^0}\right)\right)^{c_i^0/b}\right]$$

$$\text{(33)} \qquad\qquad \leq \left(\mathbf{E}\left[\exp\left(\beta \cdot \frac{b_i \tilde{x}_i(t)}{c_i^0}\right)\right]\right)^{c_i^0/b}$$

$$\text{(34)} \qquad\qquad \leq [\Psi_i(\tilde{\mathbf{s}}_i(t))]^{c_i^0/b}$$

$$\text{(35)} \qquad\qquad \leq (2e^{(1-\varepsilon/2)\beta})^{c_i^0}/b,$$

where (33) follows from Jensen's inequality applied to the concave function $x^a$, $a \leq 1$; (34) holds because $x^a$ is monotonically increasing for $a > 0$ and (35) follows from (28). From (32) and (35), we have

$$\exp\left(\frac{\beta}{b} \cdot \mathbf{E}[(\tilde{w}(t) - b)^+]\right) \leq 1 + e^{-\beta} \prod_{1 \leq i \leq m} (2e^{(1-\varepsilon/2)\beta})^{c_i^0/b}$$

$$\text{(36)} \qquad\qquad \leq 1 + e^{-\beta}(2e^{(1-\varepsilon/2)\beta})$$

$$\leq 1 + 2e^{-(\varepsilon/2)\beta},$$

where (36) follows from the bound $\sum_{i=1}^{m} c_i^0 = (1 + 4\varepsilon) \sum_{i=1}^{m} b_i \rho_i \alpha_i^\varepsilon \leq (1 + 4\varepsilon)(\frac{b}{1+4\varepsilon}) = b$. Part (i) follows by taking logarithms.

A similar argument establishes that

$$\exp\left(\beta \cdot \frac{b_i \mathbf{E}[\tilde{y}(t)]}{c_i^1}\right) \leq \mathbf{E}[e^{\beta(b_i \tilde{y}(t))/c_i^1}] \leq 2e^{(1-\varepsilon/2)\beta}.$$

Therefore,

$$(37) \quad \mathbf{E}[\tilde{y}(t)] \leq \left(\frac{\log(2)}{\beta} + 1 - \frac{\varepsilon}{2}\right)\frac{c_i^1}{b_i} \leq (1+\zeta)(1-\alpha_i^\varepsilon)\rho_i,$$

where

$$\zeta = \left(\frac{\log(2)}{\beta} + 1 - \frac{\varepsilon}{2}\right)(1+4\varepsilon) - 1.$$

Let $q_i(t)$ denote the number of class $i$ requests at time $t$ in an infinite capacity system with no admission control and let $\tilde{y}_i^0(t)$ denote the number of requests surviving from the $\tilde{y}_i(0^-)$ class $i$ requests initially loaded into system 1. Then conservation implies

$$(38) \quad q_i(t) + \tilde{y}_i^0(t) \stackrel{\mathrm{d}}{=} \tilde{x}_i(t) + \tilde{y}_i(t),$$

where $\stackrel{\mathrm{d}}{=}$ denotes equality in distribution. [Note that the surviving requests $\tilde{y}_i^0(t)$ are also counted as part of $\tilde{y}_i(t)$.] Suppose the initial load $\tilde{y}_i(0^-) = (1-\alpha_i^\varepsilon)\rho_i$, $i=1,\ldots,m$. Then

$$\frac{b_i \tilde{y}_i(0^-)}{c_i} = \frac{1}{1+4\varepsilon} \leq 1 - \frac{\varepsilon}{2} \qquad \forall\, i=1,\ldots,m;$$

that is, the hypothesis of Lemma 2 holds for all $i=1,\ldots,m$. Therefore, (37) and (38) imply that

$$(39) \quad \begin{aligned} \mathbf{E}[\tilde{x}_i(t)] &\geq \rho_i(1-\exp(-\mu_i(t))) + (1-\alpha_i^\varepsilon)\exp(-\mu_i t) - (1+\zeta)(1-\alpha_i^\varepsilon)\rho_i \\ &= \alpha_i^\varepsilon \rho_i(1-\exp(-\mu_i t)) - \zeta(1-\alpha_i^\varepsilon)\rho_i. \end{aligned}$$

Thus,

$$(40) \quad \mathbf{E}[\tilde{R}(t)] = \sum_{i=1}^m r_i \mathbf{E}[\tilde{x}_i(t)] \geq \sum_{i=1}^m \alpha_i^\varepsilon r_i \rho_i (1-\exp(-\mu_i t)) - \zeta \sum_{i=1}^m (1-\alpha_i^\varepsilon) r_i \rho_i.$$

□

Lemma 3 establishes that if $\beta \gg 1$ is admissible, the policy $\tilde{\boldsymbol{\pi}}$ does not significantly violate the capacity constraint and the associated reward rate $\mathbf{E}[\tilde{R}(t)]$ is close to the upper bound (15). The following result establishes that, on average, the policy $\tilde{\boldsymbol{\pi}}$ admits more requests than $\bar{\boldsymbol{\pi}}$.

LEMMA 4. *Fix $\varepsilon$, $\beta$, $(\mathbf{c}^0, \mathbf{c}^1)$ and the initial state $\tilde{\mathbf{y}}(0^-) = \bar{\mathbf{y}}(0^-)$. Let $\tilde{\boldsymbol{\pi}}$ and $\bar{\boldsymbol{\pi}}$ be the policies that correspond to these parameters. Then*

$$\bar{x}_i(t) \stackrel{\mathrm{d}}{\leq} \tilde{x}_i(t), \qquad \tilde{y}_i(t) \stackrel{\mathrm{d}}{\leq} \bar{y}_i(t), \qquad i = 1, \ldots, m,$$

*where $X \stackrel{\mathrm{d}}{\leq} Y$ denotes that, for all $u \geq 0$, we have $\mathbf{P}(X \geq u) \leq \mathbf{P}(Y \geq u)$.*

PROOF. The result is established by a coupling argument that employs another infeasible policy $\hat{\boldsymbol{\pi}}$ as a comparison policy.

The policies $\tilde{\boldsymbol{\pi}}$, $\bar{\boldsymbol{\pi}}$ and $\hat{\boldsymbol{\pi}}$ act on the same labeled Poisson arrival streams. Let the $k$th class $i$ arrival be labeled $(i, k)$. Let $\bar{X}_i(t)$ [resp. $\hat{X}_i(t)$] denote the set of labels of all class $i$ requests routed to system 0 by policy $\bar{\boldsymbol{\pi}}$ (resp. $\hat{\boldsymbol{\pi}}$) and still in service at time $t$.

The routing decision of the comparison policy $\hat{\boldsymbol{\pi}}$ is identical to that of the policy $\tilde{\boldsymbol{\pi}}$ unless policy $\tilde{\boldsymbol{\pi}}$ routes to system 1 (i.e., rejects) but policy $\bar{\boldsymbol{\pi}}$ routes the arrival to system 0 (i.e., accepts). Let $t$ be any time instant when this event occurs and suppose the arriving request has the label $(i, k)$. Since the policy $\hat{\boldsymbol{\pi}}$ does not face any capacity constraints, it must be that $\hat{x}_i(t^-) > \bar{x}_i(t^-)$, that is, there exists a request with label $(i, l) \in \hat{X}_i(t) \setminus \bar{X}_i(t)$. The policy $\hat{\boldsymbol{\pi}}$ admits the incoming request $(i, k)$ into system 0 by relabeling it $(i, l)$ and moves the job previously labeled $(i, l)$ to system $i$ and relabels it $(i, k)$. Clearly the policy $\hat{\boldsymbol{\pi}}$ is infeasible since the requests once routed to system 0 cannot be removed.

From the definition of the policy $\hat{\boldsymbol{\pi}}$ it is clear that $\hat{x}_i(t) \geq \bar{x}_i(t)$ and $\hat{y}_i(t) \leq \bar{y}_i(t)$. Notice that every time the policy $\hat{\boldsymbol{\pi}}$ removes a request before completion, the remaining service duration is $\exp(\mu_i)$, that is, the service duration of the request that replaces the removed request is, in distribution, identical to the remaining service duration. Therefore, the performance of the policy $\hat{\boldsymbol{\pi}}$ is, in distribution, identical to the policy $\tilde{\boldsymbol{\pi}}$. Thus, for all $u \geq 0$, we have

$$\mathbf{P}(\tilde{x}_i(t) \geq u) = \mathbf{P}(\hat{x}_i(t) \geq u) \geq \mathbf{P}(\bar{x}_i(t) \geq u),$$
$$\mathbf{P}(\tilde{y}_i(t) \geq u) = \mathbf{P}(\hat{y}_i(t) \geq u) \leq \mathbf{P}(\bar{y}_i(t) \geq u). \qquad \square$$

Let $\xi_i(t)$ [resp. $\eta_i(t)$] denote the number of class $i$ requests in system 1 at time $t$ that were rejected by the penalty function (resp. the capacity constraint). The expected value $\mathbf{E}[\xi_i(t)]$ is bounded as follows.

$$\mathbf{E}[\xi_i(t)] = \int_0^t \lambda_i \mathbf{P}\left(\frac{\partial \Psi_i(\bar{\mathbf{s}}_i(u))}{\partial x_i} > \frac{\partial \Psi_i(\bar{\mathbf{s}}_i(u))}{\partial y_i}\right) e^{-\mu(t-u)} \, du$$
$$= \int_0^t \lambda_i \mathbf{P}\left(\frac{\bar{x}_i(u)}{c_i^0} - \frac{\bar{y}_i(u)}{c_i^1} > \frac{1}{\beta b_i} \log\left(\frac{c_i^0}{c_i^1}\right)\right) e^{-\lambda(t-u)} \, du$$

$$
\begin{align}
(41) \qquad &\leq \int_0^t \lambda_i \mathbf{P}\left(\frac{\tilde{x}_i(u)}{c_i^0} - \frac{\tilde{y}_i(u)}{c_i^1} > \frac{1}{\beta b_i}\log\left(\frac{c_i^0}{c_i^1}\right)\right) e^{-\lambda(t-u)}\, du \\
&= \int_0^t \lambda_i \mathbf{P}\left(\frac{\partial \Psi_i(\tilde{\mathbf{s}}_i(u))}{\partial x_i} > \frac{\partial \Psi_i(\tilde{\mathbf{s}}_i(u))}{\partial y_i}\right) e^{-\mu(t-u)}\, du \\
(42) \qquad &= \mathbf{E}[\tilde{y}_i(t)],
\end{align}
$$

where (41) follows from $(\bar{x}_i(u))/c_i^0 - (\bar{y}_i(u))/c_i^1 \overset{d}{\leq} (\tilde{x}_i(u))/c_i^0 - (\tilde{y}_i(u))/c_i^1$.

The expected value $\mathbf{E}[\eta_i(t)]$ is bounded as follows:

$$
\begin{align}
\mathbf{E}[\eta_i(t)] &\leq \int_0^t \lambda_i \mathbf{P}\left(\sum_{i=1}^m \bar{x}_i(u) \geq b - b_i\right) e^{-\mu(t-u)}\, du \\
(43) \qquad &\leq \int_0^t \lambda_i \mathbf{P}\left(\sum_{i=1}^m \tilde{x}_i(u) \geq b - b_i\right) e^{-\mu(t-u)}\, du \\
&\leq e^{-\beta(1-b_i/b)} \int_0^t \lambda_i \mathbf{E}[e^{(\beta/b)\tilde{x}(u)}] \exp(-\mu_i(t-u))\, du \\
(44) \qquad &\leq 2 e^{-\beta(1-b_i/b)} e^{\beta(1-\varepsilon/2)} \int_0^t \lambda_i \exp(-\mu_i(t-u))\, du \\
(45) \qquad &\leq 2\rho_i e^{-\varepsilon/2(\beta-4)}(1 - \exp(-\mu_i t)),
\end{align}
$$

where (43) follows from Lemma 4, (44) follows from an argument similar to that in the proof of part (i) of Lemma 3 and (45) follows from the bound on $b_i$ implied by (23). From (42) and (45) it follows that

$$
\begin{align}
\mathbf{E}[\bar{x}_i(t)] &= \mathbf{E}[q_i(t)] + \mathbf{E}[y_{0,i}(t)] - \mathbf{E}[\bar{y}_i(t)] \\
&= \mathbf{E}[q_i(t)] + \mathbf{E}[y_{0,i}(t)] - (\mathbf{E}[\xi_i(t)] + \mathbf{E}[\eta_i(t)]) \\
(46) \qquad &\geq \mathbf{E}[q_i(t)] + \mathbf{E}[y_{0,i}(t)] - \mathbf{E}[\tilde{y}_i(t)] - 2\rho_i e^{-\varepsilon/2(\beta-4)}(1 - \exp(-\mu_i t)) \\
&= \mathbf{E}[\tilde{x}_i(t)] - 2\rho_i e^{-\varepsilon/2(\beta-4)}(1 - \exp(-\mu_i t)) \\
&\geq \alpha_i^\varepsilon \rho_i (1 - \exp(-\mu_i t)) - \zeta(1 - \alpha_i^\varepsilon) - 2\rho_i e^{-\varepsilon/2(\beta-4)}(1 - \exp(-\mu_i t)),
\end{align}
$$

where (46) follows from the bound (39) and $\zeta = (\frac{\log(2)}{\beta} + 1 - \frac{\varepsilon}{2})(1 + 4\varepsilon) - 1$. Thus, we have the following result.

THEOREM 2. *Suppose $\varepsilon < \frac{1}{4}$, $(\mathbf{c}^0, \mathbf{c}^1)$ are given by (21), $\beta$ satisfies (22) and the initial state $\bar{\mathbf{s}}(0^-) = (\mathbf{0}, \bar{\mathbf{y}}(0^-))$, with $\bar{y}_i(0^-) = (1 - \alpha_i^\varepsilon)\rho_i$, $i = 1, \ldots, m$. Then the reward rate $\bar{R}(t)$ of the penalty policy $\bar{\boldsymbol{\pi}}$ satisfies*

$$
(47) \qquad \mathbf{E}[\bar{R}(t)] \geq \max\left\{\sum_{i=1}^m \alpha_i^\varepsilon r_i \rho_i (1 - \exp(-\mu_i t)) - \zeta \sum_{i=1}^m (1 - \alpha_i^\varepsilon) r_i \rho_i \right.
$$

$$- 2e^{-\varepsilon/2(\beta-4)} \sum_{i=1}^{m} r_i\rho_i(1 - \exp(-\mu_i t)), 0 \bigg\},$$

where $\boldsymbol{\alpha}^\varepsilon$ is an optimal solution of the perturbed LP (20) and

$$\zeta = \left(\frac{\log(2)}{\beta} + 1 - \frac{\varepsilon}{2}\right)(1 + 4\varepsilon) - 1.$$

Let $L(t)$ denote the lower bound in (47). Then (15) and (47) imply that

(48)
$$\frac{\lim_{t\to\infty} L(t)}{R^*}$$
$$\geq \frac{\sum_{i=1}^{m} \alpha_i^\varepsilon r_i\rho_i - \zeta \sum_{i=1}^{m}(1-\alpha_i^\varepsilon)r_i\rho_i - 2e^{-\varepsilon/2(\beta-4)}\sum_{i=1}^{m} r_i\rho_i}{R^*}.$$

Recall that $(u^*, \mathbf{v}^*)$ denotes an optimal solution of dual LP (9). From the duality theory for LPs it follows that $(u^*, \mathbf{v}^*)$ is optimal for the dual of the perturbed LP (20) for all sufficiently small $\varepsilon$ [Luenberger (1984)], that is,

(49)    $\varepsilon_0 = \max\{\varepsilon: (u^*, \mathbf{v}^*)$ is optimal for the dual of $(20)\} > 0.$

Thus, for all $\varepsilon \leq \varepsilon_0$,

(50)
$$\sum_{i=1}^{m} \alpha_i^\varepsilon r_i\rho_i = \sum_{i=1}^{m} v_i^* + \frac{u^*b}{1+4\varepsilon}$$
$$= \left(\sum_{i=1}^{m} v_i^* + u^*b\right) - \frac{4\varepsilon}{1+4\varepsilon}(u^*b) \geq (1-4\varepsilon)R^*.$$

Since $\zeta \leq 8\varepsilon + \frac{2\log(2)}{\beta}$, (48) and (50) imply the following.

COROLLARY 1.   *Suppose $\varepsilon < \min\{\varepsilon_0, \frac{1}{4}\}$, where $\varepsilon_0$ is given by* (49), $(\mathbf{c}^0, \mathbf{c}^1)$ *are given by* (21), $\beta$ *satisfies* (22) *and* $\bar{y}_i(0^-) = (1-\alpha_i^\varepsilon)\rho_i$, $i = 1,\ldots, m$. *Then* $\bar{L} = \lim_{t\to\infty} L(t)$ *satisfies*

(51)   $\displaystyle \frac{\bar{L}}{R^*} \geq 1 - 12\varepsilon - \frac{2\log(2)}{\beta} - \left(2e^{\varepsilon/2(\beta-4)} + 8\varepsilon + \frac{2\log(2)}{\beta}\right)\frac{\sum_{i=1}^{m} r_i\rho_i}{R^*}.$

The term $\sum_{i=1}^{m} r_i\rho_i$ in (51) would appear, at first glance, to be large. However, recall that we had dropped from consideration all classes with $\alpha_i^* = 0$; therefore, $\sum_{i=1}^{m} r_i\rho_i = \sum_{\{i:\,\alpha_i^*>0\}} r_i\rho_i$, that is, the total incoming revenue rate of only the admitted classes.

Since $\varepsilon$ and $\beta$ cannot be chosen independently, the lower bound (51) implies that for every given load $\rho$ there is an optimal $\varepsilon^*(\rho)$ and a corresponding optimal lower bound $\bar{L}^*(\rho)$. The bound $\bar{L}^*(\rho)/R^* \to 1$ as $\rho \uparrow \infty$, that is, the penalty policy is optimal in the Halfin–Whitt limiting regime. This limiting result is further discussed in Section 3.3.

Next, we numerically compare the transient performance of the penalty policy $\bar{\pi}$ with the upper bound (15) for a three-class admission control problem defined by

$$
(52) \quad \boldsymbol{\lambda} = \begin{pmatrix} 40 \\ 80 \\ 60 \end{pmatrix}, \qquad \boldsymbol{\mu} = \begin{pmatrix} 0.5 \\ 2.0 \\ 0.3 \end{pmatrix},
$$

$$
\mathbf{r} = \begin{pmatrix} 1.00 \\ 0.25 \\ 0.75 \end{pmatrix}, \qquad \mathbf{b} = \begin{pmatrix} 0.10 \\ 0.15 \\ 0.55 \end{pmatrix}, \qquad b = 100.
$$

The optimal solution of the corresponding steady state LP (6) is $\alpha^* = [1, 1, 0.7818]^T$ and the optimal steady state reward $R^* = 207.2727$. The approximation parameter $\varepsilon$ was chosen by setting $\beta$ equal to the upper bound (23) and optimizing the bound (51) as a function of $\varepsilon$. The row marked Scale $\eta = 1$ in Table 1 displays the optimal $\varepsilon$, and the steady state and transient error of the optimized penalty policy. Since the lower bound $L(t) = 0$ for all sufficiently small $t$ [i.e., error $1 - ((L(t))/(R^*(t)))$ is 100%], we defined transient error $= \max\{(L(t))/R^* : t \geq 0.1/\mu_{\min}\}$.

These numerical computations were repeated for the scaled admission control problem defined by $\boldsymbol{\lambda}^{(k)} = k\lambda$, $\mathbf{r}^{(k)} = \frac{1}{k}\mathbf{r}$ and $\mathbf{b}^{(k)} = \frac{1}{k}\mathbf{b}$. The corresponding results are shown in the row marked Scale $\eta = k$ in Table 1.

From the numerical results, it is clear that as the load $\rho \uparrow \infty$, both the steady state and the transient error improve. Although the steady state error appears to converge to zero, the transient error appears to level off at approximately 15%. We believe that this is a consequence of the fact that the "target" $(\mathbf{c}^0, \mathbf{c}^1)$ is fixed instead of time-varying.

Table 1
*Comparison of bounds*

| Scale $\eta$ | Optimal $\varepsilon$ | Error (%) | |
|---|---|---|---|
| | | Steady state | Transient |
| 1 | 0.2500 | 51.3195 | 88.6202 |
| 2 | 0.2500 | 21.8708 | 61.7278 |
| 4 | 0.1838 | 17.1644 | 48.7918 |
| 8 | 0.1422 | 12.7112 | 39.3613 |
| 16 | 0.1100 | 9.3599 | 32.2373 |
| 32 | 0.0851 | 6.8943 | 26.9023 |
| 64 | 0.0659 | 5.1143 | 22.9311 |
| 128 | 0.0437 | 4.0341 | 19.2897 |
| 256 | 0.0338 | 2.8049 | 17.0118 |
| 512 | 0.0236 | 2.1991 | 15.2632 |
| 1024 | 0.0183 | 1.4909 | 14.1900 |

Regressing the scale $\eta$ on the steady state error $\bar{L}$, we obtain that

$$\eta = 4157.1 \bar{L}^{-2.1101}. \tag{53}$$

This power law paints quite a dismal picture: for steady state performance within 1% of the upper bound, the load $\rho = \mathcal{O}(10^4)$. Thus, the lower bound (51) suggests that the penalty policy is impractical for all but a small fraction of admission control applications. Fortunately, simulations (see Section 3.4) reassure us that the lower bound is quite weak and, in fact, the performance of the penalty is close to the upper bound even for moderate loads.

The numerical comparison of the bounds for a specific example is certainly not as conclusive and convincing as an analytical comparison. Nevertheless, we believe that the insights derived from this simple example would survive analytical scrutiny.

3.3. *Limiting regimes.* In this section, we investigate the performance of the policy $\bar{\pi}$ in the Halfin–Whitt limiting regime [Halfin and Whitt (1981)]. The regime of interest here is defined in terms of a scale parameter $n$ and the limiting regime is obtained as $n \uparrow \infty$. In the $n$th system,

$$
\begin{aligned}
\text{system capacity} \quad & b^{(n)} = b, \\
\text{class } i \text{ arrival rate} \quad & \lambda_i^{(n)} = n\lambda_i, \quad i = 1, \ldots, m, \\
\text{class } i \text{ service rate} \quad & \mu_i^{(n)} = \mu_i, \quad i = 1, \ldots, m, \\
\text{request size} \quad & b_i^{(n)} = \frac{b_i}{n}, \quad i = 1, \ldots, m, \\
\text{reward rate} \quad & r_i^{(n)} = \frac{r_i}{n}, \quad i = 1, \ldots, m.
\end{aligned}
\tag{54}
$$

Note that the service rates $\mu_i^{(n)}$ remain constant, that is, the system exhibits transient behavior even in the limit. In the regime defined by (54) the incoming workload $b_i^{(n)} \rho_i^{(n)}$ and the total reward rate $r_i^{(n)} \rho_i^{(n)}$ of each request class $i = 1, \ldots, m$ are independent of the scale parameter $n$, whereas the individual request size $b_i^{(n)}$ and reward rate $r_i^{(n)}$ scales down. An equivalent regime is one in which the request size remains constant but the system capacity $b^{(n)}$ scales up.

While it is plausible that appropriately thinning the incoming requests is a steady state optimal policy in the limit [Kelly (1991)], it is unlikely that thinning will perform well in the transient period. We show that the penalty policy $\bar{\pi}$ is able to control transient behavior without sacrificing steady state performance.

We will need some notation and preliminary results to enable us to state the main result of this section. Let $\boldsymbol{\pi}^{(n)}$ be any feasible policy for the $n$th system. Since $b_i^{(n)}\rho_i^{(n)} = b_i\rho_i$, for all $i = 1,\ldots,m$, the upper bound in (15) is still valid, that is,

$$\mathbf{E}[R^{\pi^{(n)}}(t)] \leq \min\left\{\sum_{i=1}^{m} r_i\rho_i(1 - \exp(-\mu_i t)),\right.$$
(55)
$$\left.\sum_{i=1}^{m} r_i\rho_i\alpha_i^*(1 - \exp(-\mu_i t)) + u^*b\exp(-\mu_{\min}t)\right\}.$$

Duality theory for LP [Luenberger (1984)] guarantees that

$$(56) \quad \sum_{i=1}^{m} \alpha_i^\varepsilon r_i\rho_i(1 - \exp(-\mu_i t)) \geq \sum_{i=1}^{m} \alpha_i^* r_i\rho_i(1 - \exp(-\mu_i t)) - \mathcal{O}(\varepsilon)$$

for all $\varepsilon \leq \varepsilon_0$, where $\varepsilon_0$ is given by (49). Fix $\varepsilon < \min\{\varepsilon_0, \frac{1}{4}\}$. Set $(\mathbf{c}^0, \mathbf{c}^1)$ using (21), set

$$\beta = \frac{2}{\varepsilon}\log\left(\frac{2}{\varepsilon}\right) + 4$$

and set

$$y_i(0^-) = (1 - \alpha_i^\varepsilon)\rho_i, \qquad i = 1,\ldots,m.$$

Define

$$(57) \quad n_0(\varepsilon) = \min\left\{n \geq 1 : \beta = \frac{2}{\varepsilon}\log\left(\frac{2}{\varepsilon}\right) + 4 \text{ satisfies } (23)\right\}.$$

Then, for all $n \geq n_0(\varepsilon)$, the bounds (56) and (47) imply that

$$(58) \quad L(t) \geq \sum_{i=1}^{m} r_i\rho_i\alpha_i^*(1 - e^{-\mu_i t}) - \mathcal{O}(\varepsilon).$$

Let $\bar{\mathbf{s}}^{(n)}(t) = (\bar{\mathbf{x}}^{(n)}(t), \bar{\mathbf{y}}^{(n)}(t))$ denote the state process and let $\bar{R}^{(n)}(t)$ denote the reward rate that corresponds to $\bar{\boldsymbol{\pi}}$ in the $n$th system. Then

$$\bar{x}_i^{(n)}(t) = \bar{x}^{(n)}(0^-)$$
$$+ A_i^x\left(\int_0^t \nu_{x,i}^{(n)}\left(\frac{1}{n}\bar{\mathbf{s}}^{(n)}(s)\right)ds\right) - D_i^x\left(\int_0^t \kappa_{x,i}^{(n)}\left(\frac{1}{n}\bar{\mathbf{s}}^{(n)}(s)\right)ds\right),$$
(59)
$$\bar{y}_i^{(n)}(t) = \bar{y}^{(n)}(0^-)$$
$$+ A_i^y\left(\int_0^t \nu_{y,i}^{(n)}\left(\frac{1}{n}\bar{\mathbf{s}}^{(n)}(s)\right)ds\right) - D_i^y\left(\int_0^t \kappa_{y,i}^{(n)}\left(\frac{1}{n}\bar{\mathbf{s}}^{(n)}(s)\right)ds\right),$$

where $\{(A_i^x, A_i^y, D_i^x, D_i^y) : 1 = 1, \ldots, m\}$ are independent rate 1 Poisson processes, the departure rates $(\kappa_{x,i}^{(n)}(\cdot), \kappa_{y,i}^{(n)}(\cdot))$, $i = 1, \ldots, m$, are given by

(60)
$$\kappa_{x,i}^{(n)}(\mathbf{s}) = n\mu_i x_i,$$
$$\kappa_{y,i}^{(n)}(\mathbf{s}) = n\mu_i y_i$$

and the arrival rates $(\nu_{x,i}^{(n)}(\cdot), \nu_{y,i}^{(n)}(\cdot))$, $i = 1, \ldots, m$, are given by

(61)
$$\nu_{x,i}^{(n)}(\mathbf{s}) = \begin{cases} n\lambda_i, & \frac{\partial \Psi_i}{\partial x_i} \leq \frac{\partial \Psi_i}{\partial y_i} \quad \text{and} \quad \sum_{j=1}^m b_j x_j(t) + \frac{1}{n} b_i \leq b, \\ 0, & \text{otherwise,} \end{cases}$$

$$\nu_{y,i}^{(n)}(\mathbf{s}) = \begin{cases} n\lambda_i, & \frac{\partial \Psi_i}{\partial x_i} > \frac{\partial \Psi_i}{\partial y_i} \quad \text{or} \\ & \frac{\partial \Psi_i}{\partial x_i} \leq \frac{\partial \Psi_i}{\partial y_i} \quad \text{and} \quad \sum_{j=1}^m b_j x_j(t) + \frac{1}{n} b_i > b, \\ 0, & \text{otherwise.} \end{cases}$$

Fix time $t$ and define $X_n = \bar{R}^{(n)}(t)$. Then

(62)
$$X_n = \sum_{i=1}^m r_i^{(n)} \bar{x}_i^{(n)} \leq \sum_{i=1}^m r_i^{(n)} \left( \frac{b}{b_i^{(n)}} \right) = b \left( \sum_{i=1}^m \frac{r_i}{b_i} \right).$$

From the dynamics (59) it follows that

(63)
$$\mathbf{var}(X_n) = \sum_{i=1}^m (r_i^{(n)})^2 \left[ \mathbf{var}\left( A_i^x \left( \int_0^t \nu_{x,i}^{(n)}\left(\frac{1}{n}\bar{\mathbf{s}}^{(n)}(s)\right) ds \right) \right) \right.$$
$$\left. + \mathbf{var}\left( D_i^x \left( \int_0^t \kappa_{x,i}^{(n)}\left(\frac{1}{n}\bar{\mathbf{s}}^{(n)}(s)\right) ds \right) \right) \right]$$
$$\leq \sum_{i=1}^m \frac{r_i}{n^2} \left( n\lambda t + n\mu_i \frac{b}{b_i} t \right).$$

The upper bounds (62) and (63) imply that the family of random variables $\{X_n : n \geq 1\}$ is tight and all its limit points are nonrandom.

To show that the sequences $\{X_n : n \geq 1\}$ have a limit, we need new notation. Let $X_q^p$ denote the reward rate at time $t$ when the policy $\bar{\pi}$ is employed in an admission control problem where the arrival rates $\lambda_i^{(p)} = p\lambda_i$, $i = 1, \ldots, m$, the capacity is $qb$ and the individual rewards $r_i$ are unscaled.

Then $X_n = \frac{1}{n}X_n^n$ and, for all $n \geq m$, one has the inequality

$$\mathbf{E}[X_n] = \frac{1}{n}\mathbf{E}[X_n^n] \geq \frac{1}{n}\mathbf{E}[X_m^n] \tag{64}$$

$$\geq \frac{1}{n}\mathbf{E}[X_m^m] = \frac{m}{n}\mathbf{E}[X_m]. \tag{65}$$

Intuitively, inequality (64) follows from the fact that the expected reward rate is a nondecreasing function of capacity, and (65) follows from the fact that, since no costs are incurred for rejecting customers, the expected reward is a nondecreasing function of the arrival rate. A formal proof of this statement requires a coupling argument very similar to that in Lemma 4.

Let $\gamma_i$, $i = 1, 2$, denote two distinct limit points of the sequence $\{X_n : n \geq 1\}$ and choose subsequences $X_{n_k} \to \gamma_1$ and $X_{m_k} \to \gamma_2$. From (62), we have $\mathbf{E}[X_{n_k}] \to \gamma_1$ and $\mathbf{E}[X_{m_k}] \to \gamma_2$. By possibly choosing subsequences, ensure that $m_k + \sqrt{m_k} \geq n_k \geq m_k$. Then (65) implies that $\gamma_1 \geq \gamma_2$. Since the order of the $\gamma_i$ was arbitrary, it follows that $\gamma_1 = \gamma_2$, that is, $\lim_{n \to \infty} X_n = X$, where $X$ is nonrandom. Thus, we have the following result.

THEOREM 3. *Suppose $\varepsilon < \min\{\varepsilon_0, \frac{1}{4}\}$, where $\varepsilon_0$ is given by (49), $(\mathbf{c}^0, \mathbf{c}^1)$ are given by (21), $\beta = \frac{2}{\varepsilon}\log(\frac{2}{\varepsilon}) + 4$ and $y_i(0^-) = (1 - \alpha_i^\varepsilon)\rho_i$, $i = 1, \ldots, m$. Let $\bar{R}^{(n)}(t)$ denote the reward rate of the policy $\bar{\pi}$ in the $n$th system. Then $\bar{R}^\infty(t) = \lim_{n \to \infty} \bar{R}^{(n)}(t)$ exists a.s. and is nonrandom. Moreover,*

$$\bar{R}^{(\infty)}(t) \geq \sum_{i=1}^{m} r_i \rho_i \alpha_i^*(1 - \exp(-\mu_i t)) - \mathcal{O}(\varepsilon), \tag{66}$$

*where $\boldsymbol{\alpha}^*$ is an optimal solution of the LP (6).*

Since the control is a discontinuous function of the state, we cannot assert that the process $\{R^{(n)}(t) : t \in [0, T]\}$ converges to the process $\{R^\infty(t) : t \in [0, T]\}$.

3.4. *Numerical experiments.* In this section we report the results of some preliminary simulation studies of the penalty policy. The objectives of these simulation experiments were to investigate the following:

(i) *The quality of the lower bound* (47): The numerical computations in Section 3.2 imply that $\rho_i = \mathcal{O}(\bar{L}^{-2.11})$ for the penalty policy to be able achieve a steady state error of order $\bar{L}$. If the lower bound were tight, this would imply that the penalty policy is impractical for all but a fraction of admission control applications. We compared the lower bound with simulated performance to evaluate the quality of the bound.

(ii) *Comparison with the thinning policy* [Kelly (1991)]: We compared the performance of the penalty and thinning policies in reward maximization and load balancing scenarios.

3.4.1. *Comparison with bounds.* We arbitrarily chose the following three scenarios.

SCENARIO 1.

$$\boldsymbol{\lambda} = \begin{pmatrix} 4 \\ 8 \\ 6 \end{pmatrix}, \qquad \boldsymbol{\mu} = \begin{pmatrix} 0.5 \\ 2 \\ 0.3 \end{pmatrix},$$

(67)

$$\mathbf{r} = \begin{pmatrix} 1 \\ 0.25 \\ 0.75 \end{pmatrix}, \qquad \mathbf{b} = \begin{pmatrix} 0.1 \\ 0.015 \\ 0.055 \end{pmatrix}, \qquad b = 1.$$

SCENARIO 2.

$$\boldsymbol{\lambda} = \begin{pmatrix} 4 \\ 8 \\ 6 \end{pmatrix}, \qquad \boldsymbol{\mu} = \begin{pmatrix} 1 \\ 2 \\ 0.3 \end{pmatrix},$$

(68)

$$\mathbf{r} = \begin{pmatrix} 1 \\ 0.25 \\ 0.75 \end{pmatrix}, \qquad \mathbf{b} = \begin{pmatrix} 0.01 \\ 0.015 \\ 0.055 \end{pmatrix}, \qquad b = 1.$$

SCENARIO 3.

$$\boldsymbol{\lambda} = \begin{pmatrix} 4 \\ 8 \\ 6 \\ 4 \end{pmatrix}, \qquad \boldsymbol{\mu} = \begin{pmatrix} 0.5 \\ 2 \\ 0.3 \\ 0.2 \end{pmatrix},$$

(69)

$$\mathbf{r} = \begin{pmatrix} 1 \\ 0.25 \\ 0.75 \\ 0.67 \end{pmatrix}, \qquad \mathbf{b} = \begin{pmatrix} 0.02 \\ 0.015 \\ 0.055 \\ 0.045 \end{pmatrix}, \qquad b = 1.$$

For each of the scenarios, the optimal solution $\boldsymbol{\alpha}^*$ and the maximum reward $R^*$ are determined by solving the LP (6). The approximation parameter $\varepsilon$ was set to the value that minimized the steady state error (51) and $\beta$ was set equal to the bound (23). The performance of the penalty policy was simulated over the period $[0, t_{\max} = 10/\mu_{\min}]$ and the reward rates were averaged over $p = 100$ independent simulation runs. The simulation was repeated for scaled systems ($\boldsymbol{\lambda}^{(n)} = n\boldsymbol{\lambda}$, $\boldsymbol{\mu}^{(n)} = \boldsymbol{\mu}$, $\mathbf{b}^{(n)} = \frac{1}{n}\mathbf{b}$, $\mathbf{r}^{(n)} = \frac{1}{n}\mathbf{r}$) $n = 10, 100, 1000$ (see Section 3.3 for details).

Figures 1–3 compare the simulation estimates with the upper bound (15) and the lower bound (47) for the three scenarios. In the plots, the reward rate is normalized by $R^*$ and time is in units of $1/\mu_{\min}$.

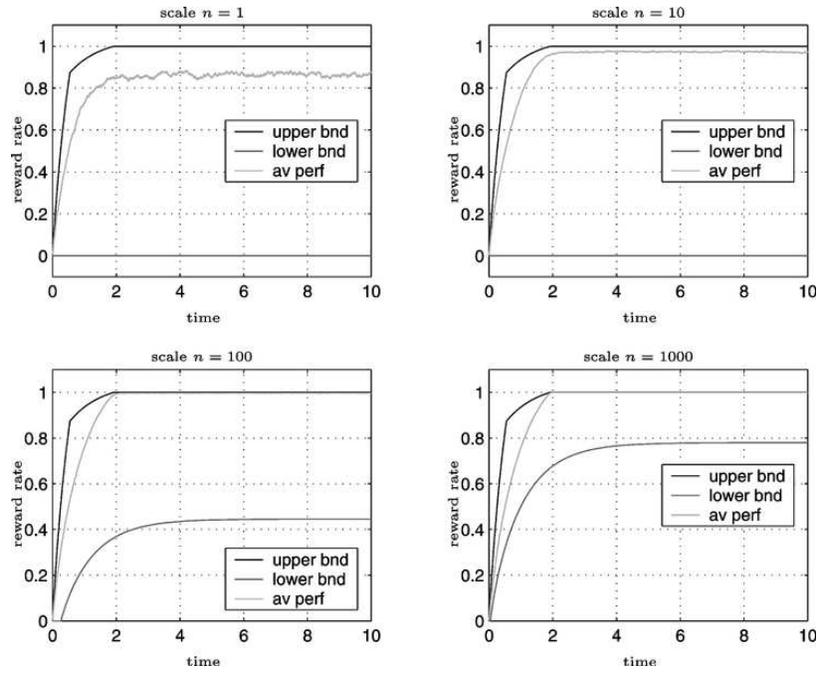

Fig. 1. *Comparison with bounds: Scenario 1.*

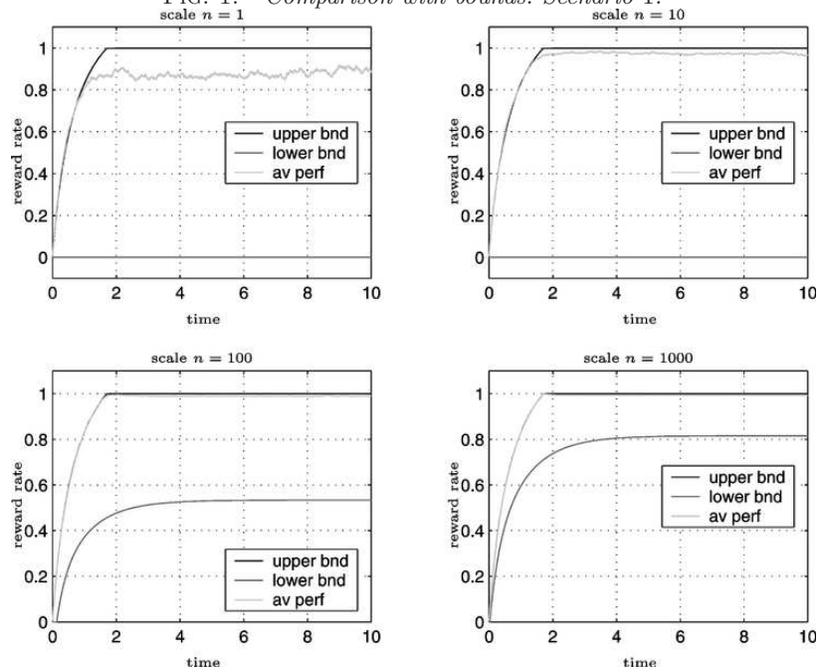

Fig. 2. *Comparison with bounds: Scenario 2.*

From the plots, it is obvious that the lower bound is quite weak, particularly so for small values of the scale parameter $n$. The performance of the penalty policy is, in fact, quite close to the upper bound. Although the transient performance of the penalty policy is significantly superior to the lower bound, it is clear that there remains a gap that needs to be bridged. Comparing the plots for different scales $n$, we see that the performance of the penalty policy is not very sensitive to the scale parameter $n$. In summary, the performance of the penalty policy, even for small loads, is remarkably good.

3.4.2. *Comparison in reward maximization scenarios.* The thinning policy is defined as follows [Kelly (1991)]. Let $\boldsymbol{\alpha}^*$ denote an optimal solution of the steady state LP (6). The thinning policy admits an arriving class $i$ request with probability $\alpha_i^*$, provided there is adequate capacity to serve the request.

Figures 4–6 plot the average performance of the penalty policy and the thinning policy as a function of the scale parameter $n$ for the three scenarios. As before, the performance was simulated over the period $[0, t_{\max} = 10/\mu_{\min}]$ and reward rates averaged over $p = 100$ independent simulation runs. In these simulation experiments both the penalty policy and the thinning policy saw the same sample path of Poisson arrivals. Also, a request accepted by both policies had the same service time in both cases.

The simulation results suggest the following conclusions. The variance of the reward rate of the thinning policy is significantly larger than the variance of the reward rate of the penalty policy. This is particularly the case for small loads. As the load increases, the steady state behavior of the thinning and penalty policies converges; however, the penalty policy remains significantly superior in the transient period.

3.4.3. *Comparison with thinning in load balancing scenarios.* The objective here is to maintain the load of the various classes close to a prescribed fraction $\mathbf{f}$, that is, class $i$ load has to be maintained close to $bf_i$, $i = 1, \ldots, m$. We considered the following two scenarios:

SCENARIO 4.

$$\text{(70)} \quad \boldsymbol{\lambda} = \begin{pmatrix} 1000 \\ 1000 \end{pmatrix}, \qquad \boldsymbol{\mu} = \begin{pmatrix} 10 \\ 10 \end{pmatrix},$$
$$\mathbf{b} = \begin{pmatrix} 1 \\ 1 \end{pmatrix}, \qquad \boldsymbol{\alpha} = \begin{pmatrix} 0.1 \\ 0.9 \end{pmatrix}, \qquad b = 100.$$

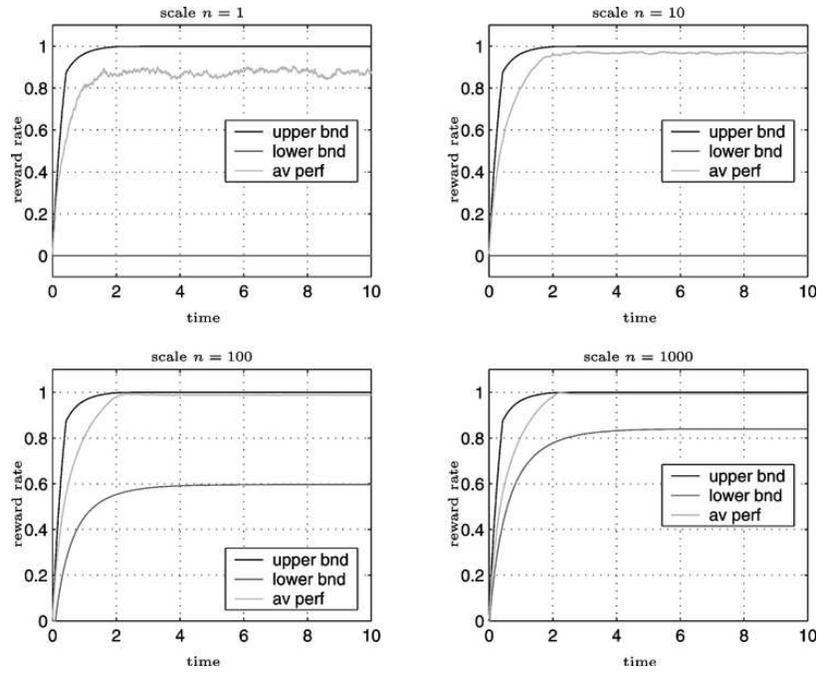

Fig. 3. *Comparison with bounds: Scenario* 3.

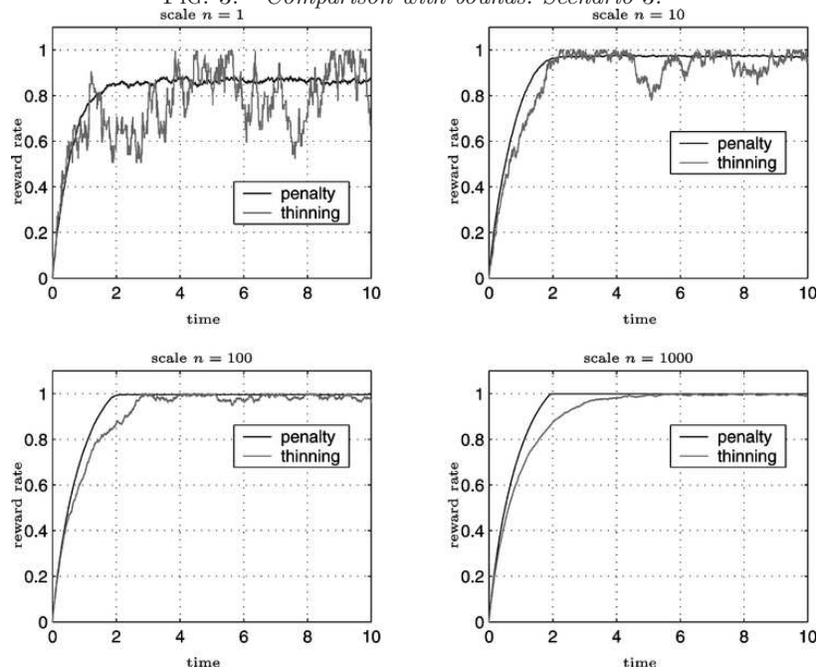

Fig. 4. *Comparison with thinning policy: Scenario* 1.

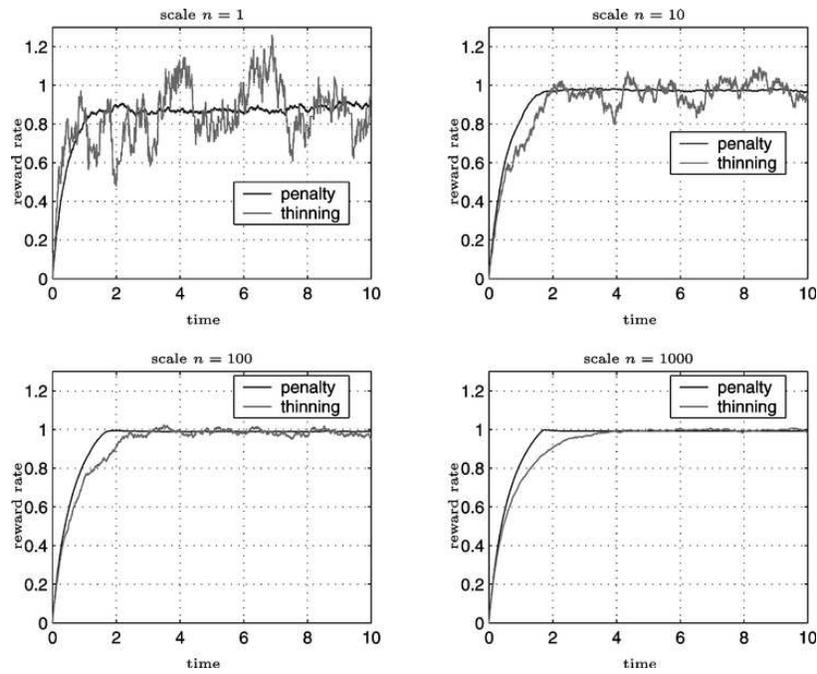

Fig. 5. *Comparison with thinning policy: Scenario 2.*

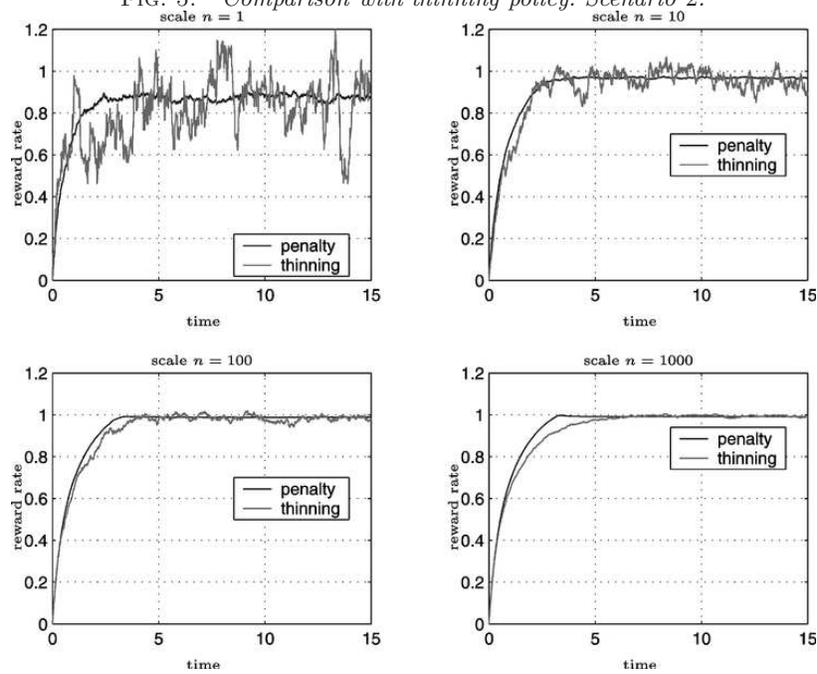

Fig. 6. *Comparison with thinning policy: Scenario 3.*

Scenario 5.

(71)
$$\boldsymbol{\lambda} = \begin{pmatrix} 100 \\ 100 \end{pmatrix}, \qquad \boldsymbol{\mu} = \begin{pmatrix} 0.1 \\ 1 \end{pmatrix},$$
$$\mathbf{b} = \begin{pmatrix} 1 \\ 1 \end{pmatrix}, \qquad \boldsymbol{\alpha} = \begin{pmatrix} 0.1 \\ 0.9 \end{pmatrix}, \qquad b = 190.$$

The two scenarios differ only in the fact that in Scenario 4, $\mu_1 = \mu_2$, whereas in Scenario 5, $\mu_2 = 10\mu_1$.

The load balancing is achieved via an appropriate admission control policy. Suppose a fraction $\alpha_i$ of all incoming class $i$ requests is admitted into the system. Then the steady state class $i$ load is $b_i \rho_i \alpha_i$. Thus, if $\alpha_i = bf_i/b_i\rho_i$, then the steady state class $i$ load will be $f_i b$. In this set of simulation experiments, we compared the performance of the thinning and penalty policies constructed from the computed admission ratio $\boldsymbol{\alpha}$.

The results for the two scenarios are shown in Figures 7 and 8. The top plot corresponds to the penalty policy and the bottom plot corresponds to the thinning policy. In both plots, the $x$-axis is time (here time is not normalized) and the $y$-axis is the fraction of the resource utilized by the requests. As before, the results are averaged over $p = 100$ iterations.

In steady state, the performance of the thinning and penalty policies is almost identical. However, the transient performance of the penalty policy is significantly superior to that of the thinning policy: In Scenario 5, where $\mu_1 \neq \mu_2$, the resource sharing that corresponds to the penalty reaches steady state levels at $t = 0.2 = 2\mu_{\min}$, whereas the resource sharing associated with the thinning policy does not reach steady state levels even by $t = 2 = 20\mu_{\min}$.

This example illustrates the target-tracking nature of the penalty policy. The policy merely tracks the target set by the capacities $(\mathbf{c}^0, \mathbf{c}^1)$. It is approximately optimal in the revenue maximization scenario because the LP sets an appropriate target to track. It could just as easily track a target set by other considerations.

3.5. *General service times.* In this section, we assume that the service duration $S_i$ has a general distribution with mean $\frac{1}{\mu_i}$, $i = 1, \ldots, m$. Let $g_i$ denote the density and let $G_i$ denote the cumulative distribution function (CDF) of the service duration $S_i$, $i = 1, \ldots, m$.

Since the steady state LP (6) and its dual (9) depend only on the mean service time $\mu_i$, they still remain the same. As before, let $R^*$ denote the optimal value, let $\boldsymbol{\alpha}^*$ denote an optimal solution of the primal LP (6) and let $(u^*, \mathbf{v}^*)$ denote an optimal solution of the dual LP (9).

Let $q_i(t)$ denote the number of active class $i$ requests at time $t$ in an infinite capacity system service time $S_i \sim g_i$ and no admission control. It is well known that [see, e.g., Wolff (1989)]

(72) $$\mathbf{E}[q_i(t)] = \rho_i(1 - \bar{G}_i^e(t)),$$

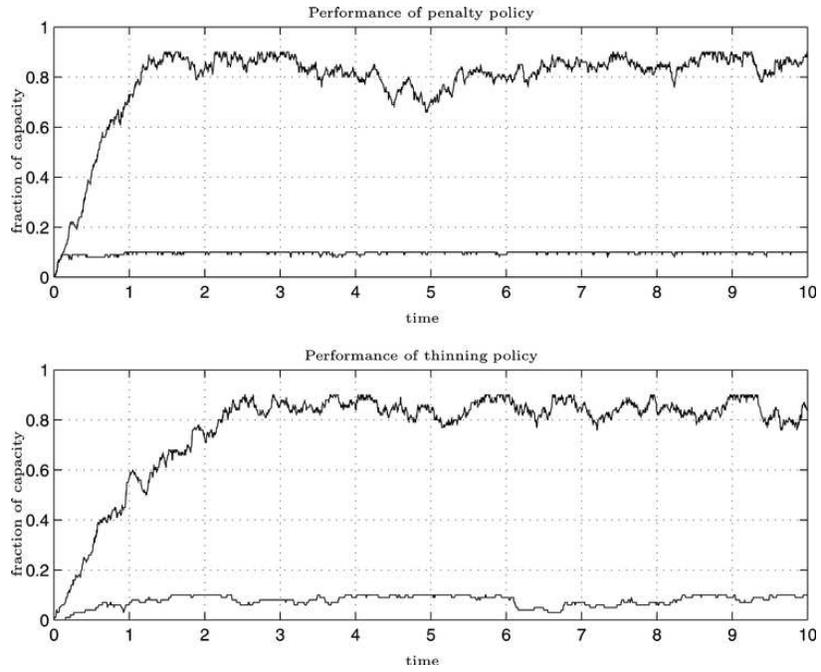

Fig. 7. *Comparison in load balancing: Scenario* 4.

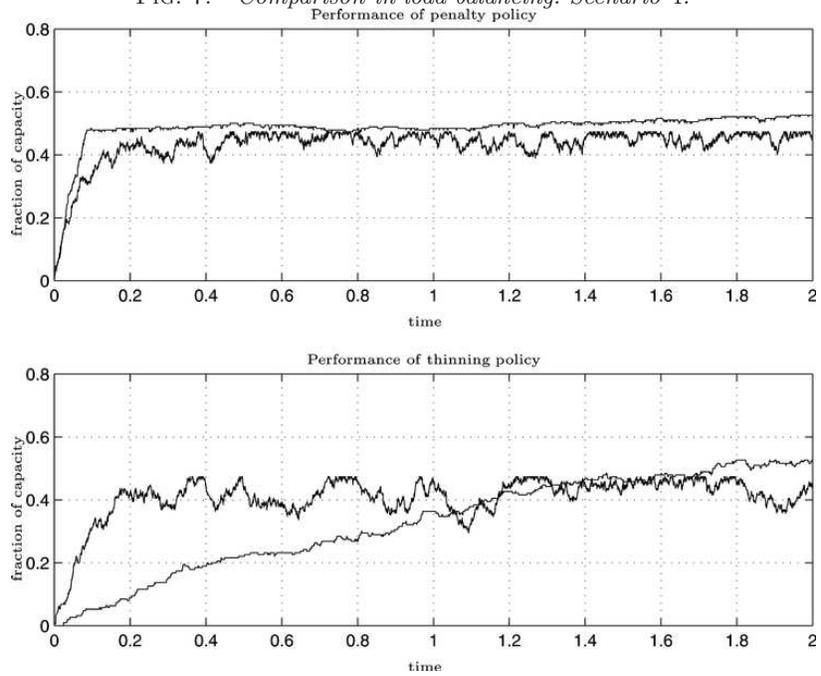

Fig. 8. *Comparison in load balancing: Scenario* 5.

where $\bar{G}_i^e(t)$ is the tail of the equilibrium CDF of the class $i$ service distribution. Thus, $\bar{G}_i^e(t)$ plays the role of the tail $\exp(-\mu_i t)$ of the exponential service time distribution. This observation leads to the following extension of Theorem 1.

THEOREM 4. *The reward rate $R^\pi(t)$ of any feasible policy $\pi$ satisfies*

$$
\begin{aligned}
\mathbf{E}[R^\pi(t)] \leq \min\bigg\{&\sum_{i=1}^m r_i \rho_i (1 - \bar{G}_i^e(t)), \\
&\sum_{i=1}^m \alpha_i^* r_i \rho_i (1 - \bar{G}_i^e(t)) + u^* b \Big(\max_{1 \leq i \leq m} \bar{G}_i^e(t)\Big)\bigg\},
\end{aligned}
\tag{73}
$$

*where $\boldsymbol{\alpha}^*$ is an optimal solution of* (6), *$(u^*, \mathbf{v}^*)$ is an optimal solution of* (9) *and $\bar{G}_i^e(\cdot)$ is the tail of the equilibrium CDF of the class $i$ service duration, $i = 1, \ldots, m$.*

Note that

$$\lim_{t \to \infty} \bigg(\sum_{i=1}^m r_i \rho_i \alpha_i^* (1 - \bar{G}_i^e(t)) + u^* b \max_{1 \leq i \leq m} \bar{G}_i^e(t)\bigg) = \sum_{i=1}^m r_i \rho_i \alpha_i^* = R^*,$$

that is, the steady state reward rate of any admissible policy is bounded above by the optimal value of the steady state LP (6).

REMARK 2. Note that in evaluating the upper bound (74), we use only the fact that the policy $\pi$ is feasible and use the bounds on the population of an $M/G/\infty$ queue [see, e.g., Wolff (1989)].

Next, we characterize the performance of the penalty policy $\bar{\boldsymbol{\pi}}$ in this model. Recall that admission decisions of the policy $\bar{\boldsymbol{\pi}}$ depend only on the load of requests of each class that have been assigned to the original system and the fictitious infinite capacity system. In particular, the policy does not keep track of the remaining service times of the requests in the system.

Let $g_i^t$ and $G_i^t$ denote, respectively, the density and the CDF of the remaining service time of a class $i$ request conditioned on the fact that it has been in service for $t$ time units. Then the tail

$$\bar{G}_i^t(s) = 1 - G_i^t(s) = \frac{G_i^e(t+s) - G_i^e(s)}{G_i^e(t)} \tag{74}$$

and, therefore,

$$g_i^t(s) = -\frac{d\bar{G}_i^t(s)}{ds} = \frac{g_i^e(s) - g_i^e(t+s)}{G_i^e(t)}. \tag{75}$$

We make the following assumption about the rate function $g_i^t(0)$.

ASSUMPTION 1. The function $g_i^t(0)$ is a decreasing function of $t$ for all $i = 1, \ldots, m$, that is, $g_i^t(0) \geq \lim_{u \to \infty} g_i^u(0) = g_i^e(0) = \mu_i$ for all $i = 1, \ldots, m$.

REMARK 3. The exponential distribution satisfies this assumption as does the heavy-tailed CDF $G(s) = (1 - (1/(1+s)^2))\mathbf{1}\{s \geq 0\}$.

Under Assumption 1, we have the following analog of Theorem 2.

THEOREM 5. *Suppose $\varepsilon < \frac{1}{4}$, $(\mathbf{c}^0, \mathbf{c}^1)$ are given by* (21), *$\beta$ satisfies* (22) *and $\bar{y}_i(0^-) = (1 - \alpha_i^\varepsilon)\rho_i$, $i = 1, \ldots, m$. Suppose also that Assumption 1 holds. Then the reward rate $\bar{R}(t)$ of the penalty policy satisfies*

$$\mathbf{E}[\bar{R}(t)] \geq \sum_{i=1}^m r_i \rho_i \alpha_i^\varepsilon (1 - \bar{G}_i^e(t)) - \sum_{i=1}^m r_i \rho_i (1 - \alpha_i^\varepsilon)(\bar{G}_i^e(t) - \bar{G}_i(t))$$

(76)
$$- \zeta \sum_{i=1}^m (1 - \alpha_i^\varepsilon) r_i \rho_i - 2 e^{-\varepsilon/2(\beta-4)} \sum_{i=1}^m r_i \rho_i (1 - e^{-\mu_i t}),$$

*where $\boldsymbol{\alpha}^\varepsilon$ is an optimal solution of the perturbed LP* (20) *and*

$$\zeta = \left(\frac{\log(2)}{\beta} + 1 - \frac{\varepsilon}{2}\right)(1 + 4\varepsilon) - 1.$$

REMARK 4. Unlike the lower bound (47), the bound (76) has a term $\sum_{i=1}^m r_i \rho_i (1 - \alpha_i^\varepsilon)(\bar{G}_i^e(t) - \bar{G}_i(t))$ that does not vanish as $\varepsilon \to 0$, that is, no matter how small the request size, this error cannot be surmounted. This term appears because the policy $\bar{\boldsymbol{\pi}}$ does not account for the remaining service times of the requests in the system.

**4. Extension to loss networks.** In this section, we extend the results of Section 3 to the network model introduced in Section 2. Recall that the stochastic system under consideration consists of a network of $s$ resources with capacity $\mathbf{b} \in \mathbf{R}_+^s$, where $b(k)$ is the capacity of resource $k = 1, \ldots, s$, and the system is initially empty. Requests for using this network of resources belong to $m$ Poisson arrival classes. Class $i$ requests have an arrival rate $\lambda_i$ and a service duration $S_i \sim \exp(\mu_i)$. They are willing to accept any capacity allocation from the set $\mathcal{B}_i = \{\mathbf{b}_{i1}, \ldots, \mathbf{b}_{il_i}\}$, $\mathbf{b}_{ij} \in \mathbf{R}_+^s$, and pay $r_i$ per unit time for the period the request is in the system.

4.1. *Upper bound on expected reward rate.* Let $\boldsymbol{\pi}$ be any feasible control policy for the stochastic problem. Let $x_{ij}^\pi(t)$ denote the number of class $i$ requests in the system at time $t$ that were assigned the capacity vector $\mathbf{b}_{ij} \in \mathcal{B}_i$.

The analog of (4) for the network setting is given by

$$
\begin{aligned}
\text{maximize} \quad & \sum_{i=1}^{m} r_i \rho_i \left( \sum_{j=1}^{l_i} \alpha_{ij} \right) \\
\text{subject to} \quad & \sum_{i=1}^{m} \rho_i \left( \sum_{j=1}^{l_i} \mathbf{b}_{ij} \alpha_{ij} \right) \leq \mathbf{b}, \\
& \sum_{j=1}^{l_i} \alpha_{ij} \leq 1 - e^{-\mu_i t}, \qquad i = 1, \ldots, m, \\
& \alpha_{ij} \geq 0, \qquad j = 1, \ldots, l_i, \ i = 1, \ldots, m.
\end{aligned}
\tag{77}
$$

Let $R^*(t)$ denote the optimal value of this LP. Taking the limit $t \to \infty$ in (77), we get the steady state LP

$$
\begin{aligned}
\text{maximize} \quad & \sum_{i=1}^{m} r_i \rho_i \left( \sum_{j=1}^{l_i} \alpha_{ij} \right) \\
\text{subject to} \quad & \sum_{i=1}^{m} \rho_i \left( \sum_{j=1}^{l_i} \mathbf{b}_{ij} \alpha_{ij} \right) \leq \mathbf{b}, \\
& \sum_{j=1}^{l_i} \alpha_{ij} \leq 1, \qquad i = 1, \ldots, m, \\
& \alpha_{ij} \geq 0, \qquad j = 1, \ldots, l_i, \ i = 1, \ldots, m.
\end{aligned}
\tag{78}
$$

Let $\boldsymbol{\alpha}^* = (\alpha_{ij}^*)_{\{j=1,\ldots,l_i, i=1,\ldots,m\}}$ denote an optimal solution and let $R^*$ denote the optimal value of (78). The dual of the steady state LP is given by

$$
\begin{aligned}
\text{minimize} \quad & \mathbf{b}^T \mathbf{u} + \mathbf{1}^T \mathbf{v} \\
\text{subject to} \quad & \rho_i r_i \leq v_i + \rho_i \mathbf{u}^T \mathbf{b}_{ij}, \qquad j = 1, \ldots, l_i, \ i = 1, \ldots, m, \\
& \mathbf{v} \geq \mathbf{0}, \qquad \mathbf{u} \geq \mathbf{0}.
\end{aligned}
\tag{79}
$$

Let $(\mathbf{u}^*, \mathbf{v}^*)$ denote an optimal solution of the dual LP (79). Then we have the following extension of Theorem 1.

THEOREM 6. *The reward rate $R^\pi(t)$ of any feasible policy $\pi$ satisfies*

$$
\begin{aligned}
\mathbf{E}[R^\pi(t)] &\leq R^*(t) \\
&\leq \min\Biggl\{ \sum_{i=1}^{m} r_i \rho_i (1 - \exp(-\mu_i t)), \\
&\qquad \sum_{i=1}^{m} r_i \rho_i \alpha_i^* (1 - \exp(-\mu_i t)) + (\mathbf{u}^*)^T \mathbf{b} \exp(-\mu_{\min} t) \Biggr\},
\end{aligned}
\tag{80}
$$

where $\alpha_i^* = \sum_{j=1}^{l_i} \alpha_{ij}^*$, $i = 1, \ldots, m$, $\boldsymbol{\alpha}^*$ is an optimal solution of steady state LP (78) and $(\mathbf{u}^*, \mathbf{v}^*)$ is an optimal solution of steady state dual LP (79).

4.2. *Penalty function and $\varepsilon$-feasible control policy.* As in the single-resource case, we drop from considerations all those capacity vectors $\mathbf{b}_{ij}$ which have the corresponding $\alpha_{ij}^* = 0$ and augment the network of systems by adding one additional fictitious infinite capacity system. The state $\mathbf{s}(t)$ of the augmented network is given by

$$\mathbf{s}(t) = (\mathbf{x}_1(t), \ldots, \mathbf{x}_m(t), \mathbf{y}(t)). \tag{81}$$

The state vector

$$\mathbf{x}_i(t) = (x_{i1}(t), \ldots, x_{il_i}(t)) \in \mathbf{Z}_+^{l_i} \tag{82}$$

describes the accepted requests, where $x_{ij}(t)$ is the number of active class $i$ request that have been assigned to $\mathbf{b}_{ij} \in \mathcal{B}_i$. The state vector $\mathbf{y}(t) = (y_1(t), \ldots, y_m(t)) \in \mathbf{Z}_+^m$, where $y_i(t)$ is the number of class $i$ requests in the fictitious system.

The penalty function $\Psi(\mathbf{s})$ is given by

$$\Psi(\mathbf{s}) = \sum_{i=1}^{m} \left[ \sum_{k=1}^{s} \underbrace{\exp\left(\beta \cdot \frac{\sum_{j=1}^{l_i} x_{ij} b_{ij}(k)}{c_{ik}^0}\right)}_{\Psi_{ik}(\mathbf{x}_i)} + \underbrace{\exp\left(\beta \cdot \frac{y_i}{c_i^1}\right)}_{\Psi_i(y_i)} \right], \tag{83}$$

where $\beta$, $(c_i^1, \{c_{ik}^0\}_{k=1}^s)$, $i = 1, \ldots, m$, are appropriately chosen constants. Let $\mathbf{s}_i = (\mathbf{x}_i, y_i)$ denote the components of the state vector that correspond to class $i$, let $\mathbf{C}^0 \in \mathbf{R}^{m \times s}$ denote the matrix $[c_{ik}^0]$ and let $\mathbf{c}^1 \in \mathbf{R}^m$ denote the vector $(c^1, \ldots, c_m^1)^T$.

The penalty policy $\bar{\boldsymbol{\pi}}$ for a loss network is defined as follows. Let $\bar{\mathbf{s}}(t) = (\bar{\mathbf{x}}_1, \ldots, \bar{\mathbf{x}}_m(t), \bar{\mathbf{y}}(t))$ denote the stochastic state process that corresponds to the policy $\bar{\boldsymbol{\pi}}$ and let $\bar{\mathbf{s}}_i = (\bar{\mathbf{x}}_i, \bar{y}_i)$. At time $t = 0^-$, the policy loads the infinite capacity system to the level $\bar{\mathbf{y}}(0^-)$. An incoming class $i$ request is conditionally accepted if

$$\min_{1 \leq j \leq l_i} \left\{ \sum_{k=1}^{s} \frac{\partial \Psi_{ik}}{\partial x_{ij}} \right\} \leq \frac{\partial \Psi_i}{\partial y_i}.$$

A conditionally accepted request is accepted and assigned to $\mathbf{b}_{ij} \in \mathcal{B}_i$ provided

$$j \in \arg \min_{1 \leq j' \leq l_i} \left\{ \sum_{k=1}^{s} \frac{\partial \Psi_{ik}}{\partial x_{ij'}} \right\}$$

and there is adequate capacity [i.e., $\sum_{i'=1}^{m} \sum_{j'=1}^{l_{i'}} \mathbf{b}_{i'j'} \bar{x}_{i'j'}(t) + \mathbf{b}_{ij} \leq \mathbf{b}$]. Otherwise the request is routed to the fictitious system and is assigned a service duration $S_i \sim \exp(\mu_i)$ that is independent of everything else.

As in the case of the single-resource problem discussed in Section 3, the capacities $(\mathbf{C}^0, \mathbf{c}^1)$ determine the following perturbed version of the steady state LP (78):

$$
\begin{aligned}
\text{maximize} \quad & \sum_{i=1}^{m} r_i \rho_i \left( \sum_{j=1}^{l_i} \alpha_{ij} \right) \\
\text{subject to} \quad & \sum_{i=1}^{m} \rho_i \left( \sum_{j=1}^{l_i} \mathbf{b}_{ij} \alpha_{ij} \right) \leq \frac{1}{1+4\varepsilon} \mathbf{b}, \\
& \sum_{j=1}^{l_i} \alpha_{ij} \leq 1, \qquad i = 1, \ldots, m, \\
& \alpha_{ij} \geq 0, \qquad j = 1, \ldots, l_i, \ i = 1, \ldots, m.
\end{aligned}
\tag{84}
$$

Let $\boldsymbol{\alpha}^\varepsilon = \{\alpha_{ij}^\varepsilon : j = 1, \ldots, l_i, i = 1, \ldots, m\}$ denote an optimal solution of (84). The capacities $(\mathbf{C}^0, \mathbf{c}^1)$ are given by

$$
\begin{aligned}
c_i^1 &= (1 + 4\varepsilon) \left( 1 - \sum_{j=1}^{l_i} \alpha_{ij}^\varepsilon \right) \rho_i, \qquad i = 1, \ldots, m, \\
c_{ik}^0 &= (1 + 4\varepsilon) \nu_k \left( \sum_{j=1}^{l_i} \alpha_{ij}^\varepsilon b_{ij}(k) \right) \rho_i, \qquad k = 1, \ldots, s, i = 1, \ldots, m,
\end{aligned}
\tag{85}
$$

where $\nu_k$ is given by

$$
\nu_k = \frac{(1/(1+4\varepsilon)) b_k}{\sum_{i=1}^{m} \sum_{j=1}^{l_i} \alpha_{ij}^\varepsilon \rho_i b_{ij}(k)}, \qquad k = 1, \ldots, s.
\tag{86}
$$

The parameter $\beta$ must satisfy the bound

$$
\beta \leq \varepsilon \min \left\{ \min_{\{(i,k) : 1 \leq i \leq m, 1 \leq k \leq s\}} \left\{ \frac{c_{ik}^0}{b_{ij}(k)} \right\}, \min_{\{i : i \in U_\varepsilon^c\}} \{c_i^1\} \right\},
\tag{87}
$$

where $U_\varepsilon^c = \{i : \sum_{j=1}^{l_i} \alpha_{ij}^\varepsilon < 1, i = 1, \ldots, m\}$.

A simple extensions of the techniques developed in Section 3 allows one to establish the following analog of Theorem 2.

THEOREM 7. *Suppose $\varepsilon < \frac{1}{4}$, $(\mathbf{C}^0, \mathbf{c}^1)$ are given by (85), $\beta$ satisfies (87) and $\bar{y}_i(0^-) = (1 - \alpha_i^\varepsilon) \rho_i$, $i = 1, \ldots, m$. Then the reward rate $\bar{R}(t)$ of the penalty policy $\bar{\pi}$ satisfies*

$$
\mathbf{E}[\bar{R}(t)] \geq \sum_{i=1}^{m} \alpha_i^\varepsilon r_i \rho_i (1 - \exp(-\mu_i t)) - \zeta \sum_{i=1}^{m} (1 - \alpha_i^\varepsilon) r_i \rho_i
\tag{88}
$$
$$
- (s+1)^2 e^{-\varepsilon/2(\beta - 4)} \sum_{i=1}^{m} r_i \rho_i (1 - \exp(-\mu_i t)),
$$

where $\alpha_i^\varepsilon = \sum_{j=1}^{l_i} \alpha_{ij}^\varepsilon$, $i=1,\ldots,m$, $\boldsymbol{\alpha}^\varepsilon$ is an optimal solution of the perturbed LP (84) and

$$\zeta = \left(\frac{\log(s+1)}{\beta} + 1 - \frac{\varepsilon}{2}\right)(1+4\varepsilon) - 1.$$

**5. Extension to general polytopic constraints.** In this section we generalize the penalty approach for admission control to a related problem of state control. Although we discuss this problem in the context of a single-resource model, the results easily extend to networks.

The stochastic model is similar to that in Section 3. Requests belong to $m$ Poisson arrival classes. Class $i$ requests have arrival rate $\lambda_i$ and service duration $S_i \sim \exp(\mu_i)$. All the requests arrive at a common infinite capacity system.

Let $\mathbf{x}(t) = (x_1(t),\ldots,x_m(t)) \in \mathbf{R}_+^m$ denote the number of requests of each class in the system at time $t$. If no control is exercized, then the expected number $\mathbf{E}[x_i(t)]$ of class $i$ requests evolves according to $\mathbf{E}[x_i(t)] = \rho_i(1 - e^{-\mu_i t})$, $i=1,\ldots,m$. Therefore, the expected steady state load is $\boldsymbol{\rho}$, where $\boldsymbol{\rho} = (\rho_1,\ldots,\rho_m)^T \in \mathbf{R}_+^m$.

Let $\mathcal{S} \subset \prod_{1 \le i \le m}[0,\rho_i]$ be a polytope defined as

(89) $$\mathcal{S} = \{\mathbf{x} : \mathbf{0} \le \mathbf{x} \le \boldsymbol{\rho}, \mathbf{D}\mathbf{x} \le \mathbf{h}\},$$

where $\mathbf{D} \in \mathbf{R}^{s \times m}$ and $\mathbf{h} \in \mathbf{R}_+^s$. We assume, without loss of generality, that $\mathbf{h} \ge \mathbf{0}$. We also assume that the interior $\mathbf{int}(\mathcal{S}) \neq \varnothing$; that is, there exists $\mathbf{x} \in \mathcal{S}$ such that $\mathbf{D}\mathbf{x} < \mathbf{d}$. In this section the objective is to construct an admission control policy that ensures that $\mathbf{x}(t) \in \mathcal{S}$ with high probability.

Define the "lifted" set

(90) $$\tilde{\mathcal{S}} = \{(\mathbf{x},\mathbf{y}) : \mathbf{0} \le \mathbf{x} \le \boldsymbol{\rho},\ \mathbf{0} \le \mathbf{y} \le \boldsymbol{\rho},\ \mathbf{D}^+\mathbf{x} + \mathbf{D}^-\mathbf{y} \le \mathbf{h} + \mathbf{D}^-\boldsymbol{\rho}\},$$

where $\mathbf{D}^+ \in \mathbf{R}^{s \times m}$ with $D_{ij}^+ = \max\{D_{ij},0\}$ and $\mathbf{D}^- \in \mathbf{R}^{s \times m}$ with $D_{ij}^- = \max\{-D_{ij},0\}$. It is clear that $\mathbf{x} \in \mathcal{S}$ implies $(\mathbf{x}, \boldsymbol{\rho} - \mathbf{x}) \in \tilde{\mathcal{S}}$. The "lifting" of the state space introduces a state space expansion that is mimicked by the control policy by adding a fictitious system to the network.

Define $(\mathbf{x}^*, \mathbf{y}^*) \in \tilde{\mathcal{S}}$ as

(91) $$(\mathbf{x}^*, \mathbf{y}^*) = \arg\min_{(\mathbf{x},\mathbf{y}) \in \tilde{\mathcal{S}}} \max_{1 \le j \le s} \left\{\frac{\mathbf{d}_j^+ \mathbf{x} + \mathbf{d}_j^- \mathbf{y}}{h_j + \mathbf{d}_j^- \boldsymbol{\rho}}\right\},$$

where $\mathbf{d}_j^+$ (resp. $\mathbf{d}_j^-$) is the $j$th row of $\mathbf{D}^+$ (resp. $\mathbf{D}^-$). Define

(92) $$\gamma^* = \max_{1 \le j \le s}\left\{\frac{\mathbf{d}_j^+ \mathbf{x}^* + \mathbf{d}_j^- \mathbf{y}^*}{h_j + \mathbf{d}_j^- \boldsymbol{\rho}}\right\} = \min_{(\mathbf{x},\mathbf{y}) \in \tilde{\mathcal{S}}} \max_{1 \le j \le s}\left\{\frac{\mathbf{d}_j^+ \mathbf{x} + \mathbf{d}_j^- \mathbf{y}}{h_j + \mathbf{d}_j^- \boldsymbol{\rho}}\right\}$$

and

$$\Psi^* = \Psi\left(\frac{(1+3\varepsilon)\mu_{\max}}{\mu_{\min}}(\mathbf{x}^*, \mathbf{y}^*)\right). \tag{93}$$

CLAIM 1. *The violation $\gamma^* < 1$.*

PROOF. By assumption, there exists $\mathbf{x} \in \mathcal{S}$ such that $\mathbf{Dx} < \mathbf{d}$, that is, $(\mathbf{d}_j^+ - \mathbf{d}_j^-)\mathbf{x} < h_j \ \forall j = 1, \ldots, s$ or, equivalently, $(\mathbf{d}_j^+ \mathbf{x} + \mathbf{d}_j^-(\boldsymbol{\rho} - \mathbf{x}))/(h_j + \mathbf{d}_j^- \boldsymbol{\rho}) < 1 \ \forall j = 1, \ldots, s$. The result follows from the fact that $\mathbf{x} \in \mathcal{S}$ implies $(\mathbf{x}, \boldsymbol{\rho} - \mathbf{x}) \in \tilde{\mathcal{S}}$. □

The quantity $\gamma^*$ is a measure of the size of the set $\tilde{\mathcal{S}}$: the smaller is the value of $\gamma^*$, the larger is the size of the set $\tilde{\mathcal{S}}$.

ASSUMPTION 2. The ratio of $\mu_{\min} = \min_{1 \leq i \leq m}\{\mu_i\}$ to $\mu_{\max} = \max_{1 \leq i \leq m}\{\mu_i\}$ is bounded below by $\gamma^*$, (i.e., $\mu_{\min}/\mu_{\max} \geq \gamma^*$).

This assumption essentially requires that the size of the target set $\tilde{\mathcal{S}}$ be comparable to the rate mismatch. If the rate mismatch is large, then the target set $\tilde{\mathcal{S}}$ cannot be too small. In particular, if all the departure rates $\mu_i$ are identical, then Assumption 2 is always satisfied. All the results in this section assume that $\mu_i$, $i = 1, \ldots, m$, satisfy Assumption 2.

As in all the previous sections, we add one fictitious system that tracks the rejected requests. Let $\mathbf{x}(t)$ [resp. $\mathbf{y}(t)$] denote the state of the original system (resp. fictitious system) at time $t$, and let $\mathbf{s}(t) = (\mathbf{x}(t), \mathbf{y}(t))$. The control policy $\tilde{\boldsymbol{\pi}}$ uses a penalty function to balance the loads of accepted and rejected customers to control the state of the system to lie in $\mathcal{S}$. The penalty function $\Psi(\mathbf{s})$ is defined as

$$\Psi(\mathbf{s}) = \sum_{j=1}^{s} \exp\left(\beta \cdot \frac{\mathbf{d}_j^+ \mathbf{x} + \mathbf{d}_j^- \mathbf{y}}{h_j + \mathbf{d}_j^- \boldsymbol{\rho}}\right), \tag{94}$$

where the multiplier $\beta$ satisfies

$$\beta \leq \varepsilon\left(\min_{1 \leq j \leq s}\{h_j + \mathbf{d}_j^- \boldsymbol{\rho}\}\right). \tag{95}$$

The policy $\tilde{\boldsymbol{\pi}}$ accepts a class $i$ request if

$$\frac{\partial \Psi}{\partial x_i} \leq \frac{\partial \Psi}{\partial y_i};$$

otherwise, the request is routed to the fictitious system and the policy $\tilde{\boldsymbol{\pi}}$ attaches to it a fictitious service time $S \sim \exp(\mu_i)$ that is independent of everything else.

We have the following analog of Lemma 2.

THEOREM 8. *Suppose $\varepsilon < \frac{1}{4}$, $\beta$ satisfies (95) and $\mathbf{E}[\Psi(\tilde{\mathbf{s}}(0))] \leq \Psi^*$, where $\Psi^*$ is given by (93). Suppose also that Assumption 2 holds. Then*

$$\mathbf{E}[\Psi(\tilde{\mathbf{s}}(t))] \leq \Psi^* \qquad \forall\, t \geq 0.$$

The following result establishes that the policy $\tilde{\boldsymbol{\pi}}$ ensures that the expected value $\mathbf{E}[\tilde{\mathbf{s}}(t)]$ of the corresponding state vector lies in an $\varepsilon$-inflation of the target set $\tilde{\mathcal{S}}$.

THEOREM 9. *Suppose $\varepsilon < \frac{1}{4}$, $\beta$ satisfies (95) and the initial state $\tilde{\mathbf{y}}(0^-)$ is chosen to ensure that $\Psi((\mathbf{0}, \tilde{\mathbf{y}}(0^-)) \leq \Psi^*$, where $\Psi^*$ is given by (93). Suppose also that Assumption 2 holds. Then, for all $t \geq 0$,*

$$(96) \quad \mathbf{d}_j \mathbf{E}[\tilde{\mathbf{x}}(t)] \leq h_j + \zeta \mathbf{d}_j^- \boldsymbol{\rho} + \mathbf{d}_j^- e^{-\mathbf{M}t}(\boldsymbol{\rho} - \tilde{\mathbf{y}}(0^-)), \qquad j = 1, \ldots, s,$$

*where*

$$\zeta = \left(\frac{\log(s)}{\beta} + 3\varepsilon\right) \quad and \quad \mathbf{M} = \mathbf{diag}(\mu_i).$$

PROOF. Repeated application of Jensen's inequality implies

$$\exp\left(\beta \max_{1 \leq j \leq s} \mathbf{E}\left\{\frac{\mathbf{d}_j^+ \tilde{\mathbf{x}}(t) + \mathbf{d}_j^- \tilde{\mathbf{y}}(t)}{h_j + \mathbf{d}_j^- \boldsymbol{\rho}}\right\}\right)$$

$$(97) \qquad \leq \exp\left(\beta \mathbf{E}\left[\max_{1 \leq j \leq s}\left\{\frac{\mathbf{d}_j^+ \tilde{\mathbf{x}}(t) + \mathbf{d}_j^- \tilde{\mathbf{y}}(t)}{h_j + \mathbf{d}_j^- \boldsymbol{\rho}}\right\}\right]\right)$$

$$\leq \mathbf{E}\left[\exp\left(\beta \max_{1 \leq j \leq s}\left\{\frac{\mathbf{d}_j^+ \tilde{\mathbf{x}}(t) + \mathbf{d}_j^- \tilde{\mathbf{y}}(t)}{h_j + \mathbf{d}_j^- \boldsymbol{\rho}}\right\}\right)\right]$$

$$\leq \mathbf{E}\Psi(\tilde{\mathbf{s}}(t))$$

$$\leq \Psi^*$$

$$\leq s e^{\beta(1+3\varepsilon)},$$

where (97) follows from the definition of $\gamma^*$ in (92). Taking logarithms, we get

$$\mathbf{d}_j^+ \mathbf{E}[\tilde{\mathbf{x}}(t)] + \mathbf{d}_j^- \mathbf{E}[\tilde{\mathbf{y}}(t)] \leq \left(\frac{\log(s)}{\beta} + 1 + 3\varepsilon\right)(h_j + \mathbf{d}_j^- \boldsymbol{\rho})$$

$$\leq (1+\zeta)(h_j + \mathbf{d}_j^- \boldsymbol{\rho}).$$

The result follows by recognizing that $\mathbf{E}[\tilde{\mathbf{x}}(t)] + \mathbf{E}[\tilde{\mathbf{y}}(t)] = (\mathbf{I} - e^{-\mathbf{M}t})\boldsymbol{\rho} + e^{-\mathbf{M}t}\tilde{\mathbf{y}}(0^-)$, where $\mathbf{M} = \mathbf{diag}(\mu_i)$. □

Theorem 9 leaves the choice of the initial loading of the fictitious systems $\tilde{\mathbf{y}}(0^-)$ open. One possible choice for $\tilde{\mathbf{y}}(0^-)$ is an optimal solution of the LP

$$
\begin{aligned}
\text{minimize} \quad & \max_{1 \leq j \leq s} \mathbf{d}_j \mathbf{M}(\boldsymbol{\rho} - \mathbf{y}) \\
\text{subject to} \quad & \mathbf{d}_j \mathbf{y} \leq (h_j + \mathbf{d}_j^- \boldsymbol{\rho})\Psi^*, \quad j = 1, \ldots, s,
\end{aligned}
\tag{98}
$$

where $\Psi^*$ is given by (93). The LP (98) minimizes the tracking error subject to the constraint that $\Psi(\mathbf{0}, \tilde{\mathbf{y}}(0^-)) \leq \Psi^*$.

Our objective in this section was to demonstrate a policy $\boldsymbol{\pi}$ that ensures that the state $\mathbf{x}^\pi(t) \in \mathcal{S}$ with high probability. Since $\mathbf{0} \leq \mathbf{E}[\mathbf{x}] \leq (\mathbf{I} - e^{-\mathbf{M}t})\boldsymbol{\rho}$, Theorem 9 states that $\mathbf{E}[\tilde{\mathbf{x}}(t)]$ lies in the set

$$(99)\, \tilde{\mathcal{S}}_\varepsilon(t) = \{\mathbf{x} : \mathbf{0} \leq \mathbf{x} \leq \boldsymbol{\rho}, \mathbf{D}\mathbf{x} \leq \mathbf{h} + \zeta(\mathbf{h} + \mathbf{D}^-\boldsymbol{\rho}) + \mathbf{D}^- e^{-\mathbf{M}t}(\boldsymbol{\rho} - \tilde{\mathbf{y}}(0^-))\},$$

where $\zeta = (\frac{\log(s)}{\beta} + 3\varepsilon)$ and $\mathbf{M} = \mathbf{diag}(\boldsymbol{\mu})$. Suppose the loads $\boldsymbol{\rho}$ are high enough such that $\beta = \frac{\log(s)}{\varepsilon}$ satisfies (95). Then $\tilde{\mathcal{S}}_\varepsilon(t)$ is an $\varepsilon$-blowup of the target set.

One might be tempted to convert this expected value result into a sample-path result by using Markov's inequality. However, such an attempt will be futile. The essential problem is that, although the policy $\tilde{\boldsymbol{\pi}}$ is able to control the accepted load, the total load of class $i$ requests is uncontrollable on a sample-path basis. Therefore, one can expect a sample-path result only if the total load is well behaved. The rest of this section investigates a limiting regime where this is the case.

Consider the limiting regime defined by (54) in Section 3.3. Choose $\varepsilon < \frac{1}{4}$ and set $\beta = \frac{1}{\varepsilon}\log(s)$. Define

$$n_0(\varepsilon) = \left\lceil \frac{\beta}{\varepsilon \min_{1 \leq j \leq s}\{h_j + \mathbf{d}_j^- \boldsymbol{\rho}\}} \right\rceil. \tag{100}$$

Then, for all $n \geq n_0(\varepsilon)$, the hypotheses of Theorems 8 and 9 are true and the corresponding bounds hold. Let $\{\tilde{\mathbf{s}}^{(n)}(t) : t \geq 0\}$ be the state process when the control policy $\tilde{\boldsymbol{\pi}}$ is employed in the $n$th system. The results in Section 3.3 imply that

$$\tilde{\mathbf{s}}^\infty(t) = \lim_{n \to \infty} \tilde{\mathbf{s}}^{(n)}(t) \tag{101}$$

exists and is nonrandom. The uniform bound on the penalty function $\Psi(\mathbf{s}^{(n)}(t)) \leq \Psi^*$ implies that the sequence $\{\tilde{\mathbf{s}}^{(n)}(t) : n \geq n_0(\varepsilon)\}$ is uniformly integrable; therefore,

$$\tilde{\mathbf{s}}^{(\infty)}(t) = \mathbf{E}[\tilde{\mathbf{s}}^{(\infty)}(t)] = \lim_{n \to \infty} \mathbf{E}[\tilde{\mathbf{s}}^{(n)}(t)], \tag{102}$$

leading to the following result.

THEOREM 10. *Fix $\varepsilon < \frac{1}{4}$, $\beta \geq \frac{1}{\varepsilon}\log(s)$ and $\tilde{\mathbf{y}}(0^-)$ such that $\Psi(\mathbf{0}, \tilde{\mathbf{y}}(0^-)) \leq \Psi^*$. Then, for all $t \geq 0$,*

$$\tilde{\mathbf{x}}^{(\infty)}(t) \in \tilde{\mathcal{S}}_\varepsilon(t)\{\mathbf{x} : \mathbf{0} \leq \mathbf{x} \leq \boldsymbol{\rho}, \mathbf{D}\mathbf{x} \leq \mathbf{h} + 4\varepsilon(\mathbf{h} + \mathbf{D}^-\boldsymbol{\rho}) \\ + \mathbf{D}^- e^{-\mathbf{M}t}(\boldsymbol{\rho} - \tilde{\mathbf{y}}(0^-))\}, \quad (103)$$

*where $\mathbf{M} = \mathbf{diag}(\mu_i)$.*

A possible choice for $\tilde{\mathbf{y}}(0^-)$ is an optimal solution of the LP (98).

**6. Concluding remarks.** In this article, we combined several disparate research ideas—mathematical programming bounds [Bertsimas, Paschalidis and Tsitsiklis (1994), Gibbens and Kelly (1995), Bertsimas and Sethuraman (2002), Bertsimas and Niño Mora (1999b) and Bertsimas and Chryssikou (1999)], state- space expansion [Kamath, Palmon and Plotkin (1998)], exponential penalty functions [Bienstock (2002)] and target tracking—to construct admission control policies. These penalty-based policies are approximately optimal when the request are sufficiently granular, that is, when the resource requested by a single request is small compared to the total capacity. The policies perform well both in the transient period as well as in steady state. The steady state performance of the penalty policy is controlled by the target supplied by a linear program, while the transient performance is controlled by a fictitious system or, equivalently, by expanding the state space. The penalty-based policies are also able to track arbitrary polyhedral target sets.

There are several issues that still remain open. From the numerical comparison of the bounds in Section 3.2 and the simulation results in Section 3.4, it is clear that in the transient period there is a gap between the performance of the control policy and the upper bound on achievable performance. This gap is probably because the capacity of the fictitious systems is too high for the transient period and, as a result, a larger fraction of the arriving requests get rejected. Thus, a possible solution would be to dynamically adapt the capacity of the fictitious systems. While this approach appears to perform well in simulation, we do not have an analytical justification for it. Also, it is unsatisfying that in the Halfin–Whitt regime we are not able to prove the convergence of the process over compact intervals (see Section 3.3). While it appears that this ought to be the case, the discontinuity in the control makes such a result hard to establish.

From the simulation results for the single-resource problem, it appears that all the benefits of the penalty policy are simply a consequence of the state space expansion that results from the addition of the fictitious systems. Further simulation experiments are planned to test this hypothesis.

In any case, state space expansion is a new technique that is worth exploring further.

In addition, there is always the issue of queuing. Building on the results developed here, Cosyn and Sigman (2004) [see also Cosyn (2003)] proposed penalty-based control policies for a finite capacity model that allows waiting and reneging. The extension to queuing networks with feedback is still open.

There are also several unresolved issues at the theoretical level. Although the exponential function allows the proofs to go through, it is not clear if it is essential to the problem. Young (1995) showed that the exponential penalty approach for packing and covering problems [see, e.g., Chapter 3 in Hochbaum (1996)] can be viewed as a derandomization approach, where, at every stage of the derandomization, one is picking a decision that minimizes a Hoeffding-type exponential bound on the probability of failure. Something similar might be at work here; that is, the admission control policy could be minimizing the worst case bound of leaving the target set. This interpretation opens the possibility that the penalty policy works because the exponential function is twisting the dynamics to make the worst sample paths most likely.

**Acknowledgment.** The authors thank the anonymous referee for helpful comments.

IEOR Department
Columbia University
500 West 120th Street, MC 4704
New York, New York 10027-6699
USA
e-mail: garud@ieor.columbia.edu
e-mail: sigman@ieor.columbia.edu